\documentclass[reqno,12pt]{amsart}
\usepackage{a4wide}
\usepackage{amsmath}
\usepackage{amssymb}
\usepackage{amsthm}
\input texdraw

\catcode`\@=11

\thinlines
\newskip\Einheit \Einheit=0.6cm
\newcount\xcoord \newcount\ycoord
\newdimen\xdim \newdimen\ydim \newdimen\PfadD@cke \newdimen\Pfadd@cke
\def\PfadDicke#1{\PfadD@cke#1 \divide\PfadD@cke by2 \Pfadd@cke\PfadD@cke \multiply\PfadD@cke by2}
\long\def\LOOP#1\REPEAT{\def\BODY{#1}\ITERATE}
\def\ITERATE{\BODY \let\next\ITERATE \else\let\next\relax\fi \next}
\let\REPEAT=\fi
\def\Punkt{\hbox{\raise-2pt\hbox to0pt{\hss\scriptsize$\bullet$\hss}}}
\def\DuennPunkt(#1,#2){\unskip
  \raise#2 \Einheit\hbox to0pt{\hskip#1 \Einheit
          \raise-2.5pt\hbox to0pt{\hss\normalsize$\bullet$\hss}\hss}}
\def\NormalPunkt(#1,#2){\unskip
  \raise#2 \Einheit\hbox to0pt{\hskip#1 \Einheit
          \raise-3pt\hbox to0pt{\hss\large$\bullet$\hss}\hss}}
\def\DickPunkt(#1,#2){\unskip
  \raise#2 \Einheit\hbox to0pt{\hskip#1 \Einheit
          \raise-4pt\hbox to0pt{\hss\Large$\bullet$\hss}\hss}}
\def\Kreis(#1,#2){\unskip
  \raise#2 \Einheit\hbox to0pt{\hskip#1 \Einheit
          \raise-4pt\hbox to0pt{\hss\Large$\circ$\hss}\hss}}
\def\Diagonale(#1,#2)#3{\unskip\leavevmode
  \xcoord#1\relax \ycoord#2\relax
      \raise\ycoord \Einheit\hbox to0pt{\hskip\xcoord \Einheit
         \unitlength\Einheit
         \line(1,1){#3}\hss}}
\def\AntiDiagonale(#1,#2)#3{\unskip\leavevmode
  \xcoord#1\relax \ycoord#2\relax 
      \raise\ycoord \Einheit\hbox to0pt{\hskip\xcoord \Einheit
         \unitlength\Einheit
         \line(1,-1){#3}\hss}}
\def\Pfad(#1,#2),#3\endPfad{\unskip\leavevmode
  \xcoord#1 \ycoord#2 \thicklines\ZeichnePfad#3\endPfad\thinlines}
\def\ZeichnePfad#1{\ifx#1\endPfad\let\next\relax
  \else\let\next\ZeichnePfad
    \ifnum#1=1
      \raise\ycoord \Einheit\hbox to0pt{\hskip\xcoord \Einheit
         \vrule height\Pfadd@cke width1 \Einheit depth\Pfadd@cke\hss}%
      \advance\xcoord by 1
    \else\ifnum#1=2
      \raise\ycoord \Einheit\hbox to0pt{\hskip\xcoord \Einheit
        \hbox{\hskip-\PfadD@cke\vrule height1 \Einheit width\PfadD@cke depth0pt}\hss}%
      \advance\ycoord by 1
    \else\ifnum#1=3
      \raise\ycoord \Einheit\hbox to0pt{\hskip\xcoord \Einheit
         \unitlength\Einheit
         \line(1,1){1}\hss}
      \advance\xcoord by 1
      \advance\ycoord by 1
    \else\ifnum#1=4
      \raise\ycoord \Einheit\hbox to0pt{\hskip\xcoord \Einheit
         \unitlength\Einheit
         \line(1,-1){1}\hss}
      \advance\xcoord by 1
      \advance\ycoord by -1
    \else\ifnum#1=5
      \advance\xcoord by -1
      \raise\ycoord \Einheit\hbox to0pt{\hskip\xcoord \Einheit
         \vrule height\Pfadd@cke width1 \Einheit depth\Pfadd@cke\hss}%
    \else\ifnum#1=6
      \advance\ycoord by -1
      \raise\ycoord \Einheit\hbox to0pt{\hskip\xcoord \Einheit
        \hbox{\hskip-\PfadD@cke\vrule height1 \Einheit width\PfadD@cke depth0pt}\hss}%
    \else\ifnum#1=7
      \advance\xcoord by -1
      \advance\ycoord by -1
      \raise\ycoord \Einheit\hbox to0pt{\hskip\xcoord \Einheit
         \unitlength\Einheit
         \line(1,1){1}\hss}
    \else\ifnum#1=8
      \advance\xcoord by -1
      \advance\ycoord by +1
      \raise\ycoord \Einheit\hbox to0pt{\hskip\xcoord \Einheit
         \unitlength\Einheit
         \line(1,-1){1}\hss}
    \fi\fi\fi\fi
    \fi\fi\fi\fi
  \fi\next}
\def\hSSchritt{\leavevmode\raise-.4pt\hbox to0pt{\hss.\hss}\hskip.2\Einheit
  \raise-.4pt\hbox to0pt{\hss.\hss}\hskip.2\Einheit
  \raise-.4pt\hbox to0pt{\hss.\hss}\hskip.2\Einheit
  \raise-.4pt\hbox to0pt{\hss.\hss}\hskip.2\Einheit
  \raise-.4pt\hbox to0pt{\hss.\hss}\hskip.2\Einheit}
\def\vSSchritt{\vbox{\baselineskip.2\Einheit\lineskiplimit0pt
\hbox{.}\hbox{.}\hbox{.}\hbox{.}\hbox{.}}}
\def\DSSchritt{\leavevmode\raise-.4pt\hbox to0pt{%
  \hbox to0pt{\hss.\hss}\hskip.2\Einheit
  \raise.2\Einheit\hbox to0pt{\hss.\hss}\hskip.2\Einheit
  \raise.4\Einheit\hbox to0pt{\hss.\hss}\hskip.2\Einheit
  \raise.6\Einheit\hbox to0pt{\hss.\hss}\hskip.2\Einheit
  \raise.8\Einheit\hbox to0pt{\hss.\hss}\hss}}
\def\dSSchritt{\leavevmode\raise-.4pt\hbox to0pt{%
  \hbox to0pt{\hss.\hss}\hskip.2\Einheit
  \raise-.2\Einheit\hbox to0pt{\hss.\hss}\hskip.2\Einheit
  \raise-.4\Einheit\hbox to0pt{\hss.\hss}\hskip.2\Einheit
  \raise-.6\Einheit\hbox to0pt{\hss.\hss}\hskip.2\Einheit
  \raise-.8\Einheit\hbox to0pt{\hss.\hss}\hss}}
\def\SPfad(#1,#2),#3\endSPfad{\unskip\leavevmode
  \xcoord#1 \ycoord#2 \ZeichneSPfad#3\endSPfad}
\def\ZeichneSPfad#1{\ifx#1\endSPfad\let\next\relax
  \else\let\next\ZeichneSPfad
    \ifnum#1=1
      \raise\ycoord \Einheit\hbox to0pt{\hskip\xcoord \Einheit
         \hSSchritt\hss}%
      \advance\xcoord by 1
    \else\ifnum#1=2
      \raise\ycoord \Einheit\hbox to0pt{\hskip\xcoord \Einheit
        \hbox{\hskip-2pt \vSSchritt}\hss}%
      \advance\ycoord by 1
    \else\ifnum#1=3
      \raise\ycoord \Einheit\hbox to0pt{\hskip\xcoord \Einheit
         \DSSchritt\hss}
      \advance\xcoord by 1
      \advance\ycoord by 1
    \else\ifnum#1=4
      \raise\ycoord \Einheit\hbox to0pt{\hskip\xcoord \Einheit
         \dSSchritt\hss}
      \advance\xcoord by 1
      \advance\ycoord by -1
    \else\ifnum#1=5
      \advance\xcoord by -1
      \raise\ycoord \Einheit\hbox to0pt{\hskip\xcoord \Einheit
         \hSSchritt\hss}%
    \else\ifnum#1=6
      \advance\ycoord by -1
      \raise\ycoord \Einheit\hbox to0pt{\hskip\xcoord \Einheit
        \hbox{\hskip-2pt \vSSchritt}\hss}%
    \else\ifnum#1=7
      \advance\xcoord by -1
      \advance\ycoord by -1
      \raise\ycoord \Einheit\hbox to0pt{\hskip\xcoord \Einheit
         \DSSchritt\hss}
    \else\ifnum#1=8
      \advance\xcoord by -1
      \advance\ycoord by 1
      \raise\ycoord \Einheit\hbox to0pt{\hskip\xcoord \Einheit
         \dSSchritt\hss}
    \fi\fi\fi\fi
    \fi\fi\fi\fi
  \fi\next}
\def\Koordinatenachsen(#1,#2){\unskip
 \hbox to0pt{\hskip-.5pt\vrule height#2 \Einheit width.5pt depth1 \Einheit}%
 \hbox to0pt{\hskip-1 \Einheit \xcoord#1 \advance\xcoord by1
    \vrule height0.25pt width\xcoord \Einheit depth0.25pt\hss}}
\def\Koordinatenachsen(#1,#2)(#3,#4){\unskip
 \hbox to0pt{\hskip-.5pt \ycoord-#4 \advance\ycoord by1
    \vrule height#2 \Einheit width.5pt depth\ycoord \Einheit}%
 \hbox to0pt{\hskip-1 \Einheit \hskip#3\Einheit 
    \xcoord#1 \advance\xcoord by1 \advance\xcoord by-#3 
    \vrule height0.25pt width\xcoord \Einheit depth0.25pt\hss}}
\def\Gitter(#1,#2){\unskip \xcoord0 \ycoord0 \leavevmode
  \LOOP\ifnum\ycoord<#2
    \loop\ifnum\xcoord<#1
      \raise\ycoord \Einheit\hbox to0pt{\hskip\xcoord \Einheit\Punkt\hss}%
      \advance\xcoord by1
    \repeat
    \xcoord0
    \advance\ycoord by1
  \REPEAT}
\def\Gitter(#1,#2)(#3,#4){\unskip \xcoord#3 \ycoord#4 \leavevmode
  \LOOP\ifnum\ycoord<#2
    \loop\ifnum\xcoord<#1
      \raise\ycoord \Einheit\hbox to0pt{\hskip\xcoord \Einheit\Punkt\hss}%
      \advance\xcoord by1
    \repeat
    \xcoord#3
    \advance\ycoord by1
  \REPEAT}
\def\Label#1#2(#3,#4){\unskip \xdim#3 \Einheit \ydim#4 \Einheit
  \def\lo{\advance\xdim by-.5 \Einheit \advance\ydim by.5 \Einheit}%
  \def\llo{\advance\xdim by-.25cm \advance\ydim by.5 \Einheit}%
  \def\loo{\advance\xdim by-.5 \Einheit \advance\ydim by.25cm}%
  \def\o{\advance\ydim by.25cm}%
  \def\ro{\advance\xdim by.5 \Einheit \advance\ydim by.5 \Einheit}%
  \def\rro{\advance\xdim by.25cm \advance\ydim by.5 \Einheit}%
  \def\roo{\advance\xdim by.5 \Einheit \advance\ydim by.25cm}%
  \def\l{\advance\xdim by-.30cm}%
  \def\r{\advance\xdim by.30cm}%
  \def\lu{\advance\xdim by-.5 \Einheit \advance\ydim by-.6 \Einheit}%
  \def\llu{\advance\xdim by-.25cm \advance\ydim by-.6 \Einheit}%
  \def\luu{\advance\xdim by-.5 \Einheit \advance\ydim by-.30cm}%
  \def\u{\advance\ydim by-.30cm}%
  \def\ru{\advance\xdim by.5 \Einheit \advance\ydim by-.6 \Einheit}%
  \def\rru{\advance\xdim by.25cm \advance\ydim by-.6 \Einheit}%
  \def\ruu{\advance\xdim by.5 \Einheit \advance\ydim by-.30cm}%
  #1\raise\ydim\hbox to0pt{\hskip\xdim
     \vbox to0pt{\vss\hbox to0pt{\hss$#2$\hss}\vss}\hss}%
}
\catcode`\@=12

\numberwithin{equation}{section}

\newcounter{saveeqn}
\newcommand{\alphaeqn}{\setcounter{saveeqn}{\value{equation}}%
\setcounter{equation}{0}%
\global\def\theequation{\mbox{\thesection.\arabic{saveeqn}\alph{equation}}}}
\newcommand{\reseteqn}{\setcounter{equation}{\value{saveeqn}}%
\global\def\theequation{\thesection.\arabic{equation}}}

\newtheorem{theorem}{Theorem}
\newtheorem{conjecture}[theorem]{Conjecture}

\theoremstyle{definition}
\newtheorem{definition}[theorem]{Definition}

\theoremstyle{remark}

\catcode`\@=11
\def\sideset#1#2#3{%
  \@mathmeasure\z@\displaystyle{#3}%
  \global\setbox\@ne\vbox to\ht\z@{}\dp\@ne\dp\z@
  \setbox\tw@\box\@ne
  \@mathmeasure4\displaystyle{\copy\tw@#1}%
  \@mathmeasure6\displaystyle{#3{#2}}%
  \dimen@-\wd6 \advance\dimen@\wd4 \advance\dimen@\wd\z@
  \hbox to\dimen@{}\mathop{\kern-\dimen@\box4\box6}%
}
\catcode`\@=12

\def\ringerl(#1 #2){\move(#1 #2)\fcir f:0 r:.1}
\def\rdSchritt{\bsegment
  \lpatt(.05 .13)
  \rlvec(0.288675 -.5) 
  \savepos(0.288675 -.5 )(*ex *ey)
        \esegment
  \move(*ex *ey)
        }
\def\ldSchritt{\bsegment
  \lpatt(.05 .13)
  \rlvec(-0.288675 -.5) 
  \savepos(-0.288675 -.5 )(*ex *ey)
        \esegment
  \move(*ex *ey)
        }
\def\ruSchritt{\bsegment
  \lpatt(.05 .13)
  \rlvec(0.288675 .5) 
  \savepos(0.288675 .5)(*ex *ey)
        \esegment
  \move(*ex *ey)
        }
\def\luSchritt{\bsegment
  \lpatt(.05 .13)
  \rlvec(-0.288675 .5) 
  \savepos(-0.288675 .5)(*ex *ey)
        \esegment
  \move(*ex *ey)
        }
\def\hhSchritt{\bsegment
  \lpatt(.05 .13)
  \rlvec(0.57735 0) 
  \savepos(0.57735 0)(*ex *ey)
        \esegment
  \move(*ex *ey)
        }
\def\hhhSchritt{\bsegment
  \lpatt(.05 .13)
  \rlvec(-0.57735 0) 
  \savepos(-0.57735 0)(*ex *ey)
        \esegment
  \move(*ex *ey)
        }
\def\RDSchritt{\bsegment
  \rlvec(0.288675 -.5) 
  \savepos(0.288675 -.5 )(*ex *ey)
        \esegment
  \move(*ex *ey)
        }
\def\LDSchritt{\bsegment
  \rlvec(-0.288675 -.5) 
  \savepos(-0.288675 -.5 )(*ex *ey)
        \esegment
  \move(*ex *ey)
        }
\def\RUSchritt{\bsegment
  \rlvec(0.288675 .5) 
  \savepos(0.288675 .5)(*ex *ey)
        \esegment
  \move(*ex *ey)
        }
\def\LUSchritt{\bsegment
  \rlvec(-0.288675 .5) 
  \savepos(-0.288675 .5)(*ex *ey)
        \esegment
  \move(*ex *ey)
        }
\def\HHSchritt{\bsegment
  \rlvec(0.57735 0) 
  \savepos(0.57735 0)(*ex *ey)
        \esegment
  \move(*ex *ey)
        }
\def\hdSchritt{\bsegment
  \lpatt(.05 .13)
  \rlvec(0.866025403784439 .5) 
  \savepos(0.866025403784439 .5)(*ex *ey)
        \esegment
  \move(*ex *ey)
        }
\def\vdSchritt{\bsegment
  \lpatt(.05 .13)
  \rlvec(0 1) 
  \savepos(0 1)(*ex *ey)
        \esegment
  \move(*ex *ey)
        }
\def\ldreieck{\bsegment
  \rlvec(0.866025403784439 .5) \rlvec(0 -1)
  \rlvec(-0.866025403784439 .5)
  \savepos(0.866025403784439 -.5)(*ex *ey)
        \esegment
  \move(*ex *ey)
        }
\def\rdreieck{\bsegment
  \rlvec(0.866025403784439 -.5) \rlvec(-0.866025403784439 -.5)  \rlvec(0 1)
  \savepos(0 -1)(*ex *ey)
        \esegment
  \move(*ex *ey)
        }
\def\rhombus{\bsegment
  \rlvec(0.866025403784439 .5) \rlvec(0.866025403784439 -.5)
  \rlvec(-0.866025403784439 -.5)  \rlvec(0 1)
  \rmove(0 -1)  \rlvec(-0.866025403784439 .5)
  \savepos(0.866025403784439 -.5)(*ex *ey)
        \esegment
  \move(*ex *ey)
        }
\def\RhombusA{\bsegment
  \rlvec(0.866025403784439 .5) \rlvec(0.866025403784439 -.5)
  \rlvec(-0.866025403784439 -.5) \rlvec(-0.866025403784439 .5)
  \savepos(0.866025403784439 -.5)(*ex *ey)
        \esegment
  \move(*ex *ey)
        }
\def\RhombusB{\bsegment
  \rlvec(0.866025403784439 .5) \rlvec(0 -1)
  \rlvec(-0.866025403784439 -.5) \rlvec(0 1)
  \savepos(0 -1)(*ex *ey)
        \esegment
  \move(*ex *ey)
        }
\def\RhombusC{\bsegment
  \rlvec(0.866025403784439 -.5) \rlvec(0 -1)
  \rlvec(-0.866025403784439 .5) \rlvec(0 1)
  \savepos(0.866025403784439 -.5)(*ex *ey)
        \esegment
  \move(*ex *ey)
        }

\catcode`\@=11
\def\iddots{\mathinner{\mkern1mu\raise\p@\hbox{.}\mkern2mu
    \raise4\p@\hbox{.}\mkern2mu\raise7\p@\vbox{\kern7\p@\hbox{.}}\mkern1mu}}
\catcode`\@=12

\def\de{\delta}

\def\la{\lambda}
\def\rh{\rho}

\def\({\left(}
\def\){\right)}
\def\[{\left[}
\def\]{\right]}

\begin{document}

\title{Plane partitions in the work of Richard Stanley and his school}

\author{C. Krattenthaler}

\address{Fakult\"at f\"ur Mathematik, Universit\"at Wien,
Oskar-Morgenstern-Platz~1, A-1090 Vienna, Austria.
WWW: \tt http://www.mat.univie.ac.at/\~{}kratt.}




\begin{abstract}
These notes provide a survey of the theory of plane partitions,
seen through the glasses of the work of Richard Stanley and his
school.
\end{abstract}

\maketitle

\section{Introduction}
Plane partitions were introduced to (combinatorial) mathematics
by Major Percy Ale\-xan\-der MacMahon \cite{MacMAB} around 1900.
What he had in mind was a planar analogue of 
a(n integer) partition.\footnote{In particular, a plane partition
is {\it not\/} a partition of the plane.}

In order to explain this, let us start with the definition of
a(n integer) partition.
A {\it partition of a positive integer $n$} is a way to represent $n$
as a sum of positive integers, where the order of the summands does
not play any role. So, a partition of $n$ is
\begin{equation} \label{eq:1}
n=\la_1+\la_2+\dots+\la_k,
\end{equation}
for some $k$, where all summands $\la_i$ are positive integers and,
since the order of summands is irrelevant, we may assume without loss
of generality that $\la_1\ge\la_2\ge\dots\ge\la_k>0$. For example, there
are $7$ partitions of $5$, namely
$$
5=4+1=3+2=3+1+1=2+2+1=2+1+1+1=1+1+1+1+1.
$$
If we want to represent a partition as in \eqref{eq:1} 
{\it very} succinctly, then we just write
$$
\la_1\,\la_2\,\dots\,\la_k.
$$
A plane partition is a planar analogue of this.
A {\it plane partition of a positive integer $n$} 
is a planar array $\pi$ of non-negative integers of the form
\begin{equation} \label{eq:PP}
\begin{matrix}
\pi_{1,1}&\multicolumn{4}{c}{\dotfill} &\pi_{1,\la_1}\\
\pi_{2,1}&\multicolumn{3}{c}{\dotfill}&
\pi_{2,\la_2}\\
\vdots&&&\iddots&\\
\pi_{r,1}&\dots &\pi_{r,\la_r}
 \end{matrix}
\end{equation}
such that entries along rows {\it and\/} along 
columns are weakly decreasing, i.e.,
\begin{equation} \label{e2.7}
\pi_{i,j}\ge \pi_{i,j+1}\quad \text {and}\quad \pi_{i,j}\ge
\pi_{i+1,j},
\end{equation}
and such that the sum $\sum \pi_{i,j}$ of all entries $\pi_{i,j}$
equals $n$. Here, the sequence of row-lengths, $(\la_1,\la_2,\dots,\la_r)$,
is assumed to form a partition, i.e., $\la_1\ge \la_2\ge \dots\ge \la_r>0$,
and this partition is called the {\it shape} of the plane partition~$\pi$.
The individual entries $\pi_{i,j}$ are called {\it parts} of $\pi$.
The sum $\sum \pi_{i,j}$ of all parts of the plane partition $\pi$ ---
so-to-speak the ``size" of $\pi$ ---
i.e., the number that $\pi$ partitions, will be denoted by $\vert\pi\vert$.
For example, Figure~\ref{fig:PP} shows a plane partition of $24$. 
That is, if $\pi_0$ denotes this plane partition, then
$\vert\pi_0\vert=24$, and its shape is $(4,3,2)$.

\begin{figure}[h]
$$\begin{matrix}
5&3&3&2\\
5&1&1\\
3&1
\end{matrix}
$$
\label{fig:PP}
\caption{A plane partition}
\end{figure}

\begin{figure}[h]
\centertexdraw{
  \drawdim truecm  \linewd.02
\move(8 0)
\bsegment
  \drawdim truecm  \linewd.02
  \linewd.12
  \move(0 0)
  \RhombusA \RhombusB \RhombusB 
  \RhombusA \RhombusA \RhombusB \RhombusA \RhombusB \RhombusB
  \move (-0.866025403784439 -.5)
  \RhombusA \RhombusB \RhombusB \RhombusB \RhombusB
  \RhombusA \RhombusA \RhombusB 
  \move (-1.732050807568877 -3)
  \RhombusA \RhombusB \RhombusB \RhombusA
  \RhombusB 
  \move (-0.866025403784439 -.5)
  \RhombusC
  \move (-0.866025403784439 -1.5)
  \RhombusC
  \move (0.866025403784439 -2.5)
  \RhombusC \RhombusC 
  \move (0.866025403784439 -3.5)
  \RhombusC \RhombusC \RhombusC 
  \move (2.598076211353316 -5.5)
  \RhombusC 
  \move (0.866025403784439 -5.5)
  \RhombusC 
  \move (-1.732050807568877 -3)
  \RhombusC 
  \move (-1.732050807568877 -4)
  \RhombusC 
  \move (-1.732050807568877 -5)
  \RhombusC \RhombusC 
\linewd.05
  \move (0.866025403784439 .5)
  \avec (0.866025403784439 2.5)
  \move (-1.732050807568877 -6)
  \avec (-3.4641 -7)
  \move (4.330125 -6.5)
  \avec (6.062175 -7.5)
\esegment}
\caption{A pile of cubes}
\label{fig:PP2}
\end{figure}

There is a nice way to represent a plane partition as a 
three-dimensional object: this is done by replacing each
part $k$ of the plane partition by a stack of $k$ unit
cubes. Thus we obtain a pile of unit cubes. The pile of cubes
corresponding to the plane partition in Figure~\ref{fig:PP}
is shown in Figure~\ref{fig:PP2}. It is placed into a coordinate
system (the coordinate axes being indicated by arrows in the figure)
so that its ``back corner" resides in the origin of the coordinate
system and it aligns with the coordinate axes as shown in the figure.
These piles are not
arbitrary piles; clearly they inherit the monotonicity condition
on the parts of the original plane partition. 

If we forget that
Figure~\ref{fig:PP2} is meant to show a 3-dimensional object, 
but rather consider it as graphic object in the plane
(which it really is in this printed form),
then we realise that this object consists entirely of
rhombi with side length~1 (say) and angles of $60^\circ$ and
$120^\circ$, respectively, which fit together perfectly, in the
sense that there are no ``holes" in the figure which are not such rhombi.
By adding a few more rhombi, we may enlarge this figure to
a hexagon, see Figure~\ref{fig:PP3}.b. (At this point, the thin lines
should be ignored). This hexagon has the property that opposite sides
have the same length, and its angles are all $120^\circ$.
The upshot of all this is that plane partitions
with $a$ rows and $b$ columns, and with 
parts bounded above by $c$ (recall that parts are non-negative integers) 
are in bijection with tilings of a hexagon
with side lengths $a,b,c,a,b,c$ into unit rhombi, cf.\ again
Figures~\ref{fig:PP}--\ref{fig:PP3}, with Figure~\ref{fig:PP3}.a
showing the assignments of side lengths of the hexagon. 
This is an important point of view which has only been observed
relatively recently \cite{DaToAA}. In particular, it was not known at
the time when Richard Stanley studied plane partitions. 
Nevertheless, I shall use this point of view in the exposition
here, since it facilitates many explanations and arguments
enormously. I would even say that, had it been known earlier, much of the
plane partition literature could have been presented much more
elegantly \dots 

\begin{figure}[h]
\centertexdraw{
  \drawdim truecm  \linewd.02
  \rhombus \rhombus \rhombus \rhombus \ldreieck
  \move (-0.866025403784439 -.5)
  \rhombus \rhombus \rhombus \rhombus \rhombus \ldreieck
  \move (-1.732050807568877 -1)
  \rhombus \rhombus \rhombus \rhombus \rhombus \rhombus \ldreieck
  \move (-1.732050807568877 -1)
  \rdreieck
  \rhombus \rhombus \rhombus \rhombus \rhombus \rhombus \ldreieck
  \move (-1.732050807568877 -2)
  \rdreieck
  \rhombus \rhombus \rhombus \rhombus \rhombus \rhombus \ldreieck
  \move (-1.732050807568877 -3)
  \rdreieck
  \rhombus \rhombus \rhombus \rhombus \rhombus \rhombus 
  \move (-1.732050807568877 -4)
  \rdreieck
  \rhombus \rhombus \rhombus \rhombus \rhombus 
  \move (-1.732050807568877 -5)
  \rdreieck
  \rhombus \rhombus \rhombus \rhombus 
\move(8 0)
\bsegment
  \drawdim truecm  \linewd.02
  \rhombus \rhombus \rhombus \rhombus \ldreieck
  \move (-0.866025403784439 -.5)
  \rhombus \rhombus \rhombus \rhombus \rhombus \ldreieck
  \move (-1.732050807568877 -1)
  \rhombus \rhombus \rhombus \rhombus \rhombus \rhombus \ldreieck
  \move (-1.732050807568877 -1)
  \rdreieck
  \rhombus \rhombus \rhombus \rhombus \rhombus \rhombus \ldreieck
  \move (-1.732050807568877 -2)
  \rdreieck
  \rhombus \rhombus \rhombus \rhombus \rhombus \rhombus \ldreieck
  \move (-1.732050807568877 -3)
  \rdreieck
  \rhombus \rhombus \rhombus \rhombus \rhombus \rhombus 
  \move (-1.732050807568877 -4)
  \rdreieck
  \rhombus \rhombus \rhombus \rhombus \rhombus 
  \move (-1.732050807568877 -5)
  \rdreieck
  \rhombus \rhombus \rhombus \rhombus 
  \linewd.12
  \move(0 0)
  \RhombusA \RhombusB \RhombusB 
  \RhombusA \RhombusA \RhombusB \RhombusA \RhombusB \RhombusB
  \move (-0.866025403784439 -.5)
  \RhombusA \RhombusB \RhombusB \RhombusB \RhombusB
  \RhombusA \RhombusA \RhombusB \RhombusA 
  \move (-1.732050807568877 -1)
  \RhombusB \RhombusB \RhombusA \RhombusB \RhombusB \RhombusA
  \RhombusB \RhombusA \RhombusA 
  \move (1.732050807568877 0)
  \RhombusC \RhombusC \RhombusC 
  \move (1.732050807568877 -1)
  \RhombusC \RhombusC \RhombusC 
  \move (3.464101615137755 -3)
  \RhombusC 
  \move (-0.866025403784439 -.5)
  \RhombusC
  \move (-0.866025403784439 -1.5)
  \RhombusC
  \move (0.866025403784439 -2.5)
  \RhombusC \RhombusC 
  \move (0.866025403784439 -3.5)
  \RhombusC \RhombusC \RhombusC 
  \move (2.598076211353316 -5.5)
  \RhombusC 
  \move (0.866025403784439 -5.5)
  \RhombusC 
  \move (-1.732050807568877 -3)
  \RhombusC 
  \move (-1.732050807568877 -4)
  \RhombusC 
  \move (-1.732050807568877 -5)
  \RhombusC \RhombusC 
\esegment
\htext (-1.5 -9){\small a. A hexagon with sides $a,b,c,a,b,c$,}
\htext (-1.5 -9.5){\small \hphantom{a. }where $a=3$, $b=4$, $c=5$}
\htext (6.8 -9){\small b. A rhombus tiling of a hexagon}
\htext (6.8 -9.5){\small \hphantom{b. }with sides $a,b,c,a,b,c$}
\rtext td:0 (4.3 -4.1){$\sideset {} c 
    {\left.\vbox{\vskip2.6cm}\right\}}$}
\rtext td:60 (2.6 -.55){$\sideset {} {} 
    {\left.\vbox{\vskip2cm}\right\}}$}
\rtext td:120 (-.34 -.2){$\sideset {}  {}  
    {\left.\vbox{\vskip1.7cm}\right\}}$}
\rtext td:0 (-2.4 -3.6){$\sideset {c}  {} 
    {\left\{\vbox{\vskip2.6cm}\right.}$}
\rtext td:240 (-0.1 -6.9){$\sideset {}  {}  
    {\left.\vbox{\vskip2cm}\right\}}$}
\rtext td:300 (2.9 -7.3){$\sideset {}  {}  
    {\left.\vbox{\vskip1.7cm}\right\}}$}
\htext (-.9 0.2){$a$}
\htext (2.8 -.1){$b$}
\htext (3.2 -7.9){$a$}
\htext (-0.4 -7.65){$b$}
}
\caption{}
\label{fig:PP3}
\end{figure}

\section{Percy Alexander MacMahon}
\label{sec:MM}
The problem that MacMahon posed to himself was:
\begin{gather*}
\text{\it Given positive integers $a,b,c$,
compute the generating function $\sum_{\pi}q^{\vert\pi\vert}$,}\\
\text{\it where the
sum is over all plane partitions $\pi$ contained in an
$a\times b\times c$ box!}
\end{gather*}
As explained in the previous section,
``contained in an $a\times b\times c$ box" has to be understood
in the sense of the ``pile of cubes"-interpretation (see
Figure~\ref{fig:PP2}) of plane partitions. In terms of rhombus
tilings, we are considering rhombus tilings of a hexagon with side
lengths $a,b,c,a,b,c$ (see Figure~\ref{fig:PP3}), 
while in terms of the original
definition \eqref{eq:PP}/\eqref{e2.7}, we are considering
plane partitions of shape $(b,b,\dots,b)$ (with $a$ occurrences of
$b$) with parts at most $c$.

Why were plane partitions so fascinating for MacMahon, and for
legions of followers?  From his writings, it is clear that MacMahon
did not have any external motivation to consider these objects, nor
did he have any second thoughts. For him it was obvious that these
plane partitions are very natural, as two-dimensional analogues of
(linear) partitions (for which at the time already a well established
theory was available), and as such of intrinsic interest. Moreover,
this intuition was ``confirmed" by the extremely elegant product
formula in Theorem~\ref{thm:MM} below. He himself --- conjecturally
--- found another intriguing product formula for so-called
``symmetric" plane partitions contained in a given box (see \eqref{eq:C2a}).
Later many more such formulae were found (again, first conjecturally,
and some of them still
quite mysterious); see Section~\ref{sec:10} below. Moreover, 
over time it turned out that plane partitions (and rhombus tilings)
are related to many other areas of mathematics, most notably to
the theory of symmetric functions and 
representation theory of classical groups
(as Richard Stanley pointed out; see Sections~\ref{sec:methods}
and \ref{sec:SCPP}),
representation theory of quantum groups (cf.\ \cite{KupeAF}), 
enumeration of integer points in polytopes and commutative algebra
(cf.\ \cite{BrZaAA,PakIAZ,CGJLAA} and the references therein),
enumeration of matchings in graphs (cf.\ \cite{KupeAA,KupeAG}),
and to statistical physics (cf.\ \cite{GuOVAA,KrGVAA}).
In brief, the theory of plane partitions offered challenging
conjectures (most of them solved now, but see Section~\ref{sec:endofstory}),
connections to many other areas, and therefore many different
views and approaches to attack problems are possible, making
this a very rich research field.

Returning to MacMahon, the surprisingly simple 
formula \cite[Sec.~429]{MacMAA} he found for the generating
function for all plane partitions contained in a given box
was the following.

\begin{theorem} \label{thm:MM}
The generating function $\sum_\pi q^{\vert\pi\vert}$ for all 
plane partitions $\pi$ contained in an $a\times b\times c$ box
is given by\footnote{MacMahon did not write the result in this form,
for which there are many ways to express it. The particular product
\eqref{eq:M2} is due to Macdonald \cite[Eq.~(2) on p.~81]{MacdAC}.}
\begin{equation} \label{eq:M2}
\prod _{i=1} ^{a}
\prod _{j=1} ^{b}
\prod _{k=1} ^{c}
\frac {1-q^{i+j+k-1}} {1-q^{i+j+k-2}}.
\end{equation}
\end{theorem}

We note that, if we let $a,b,c\to\infty$, then we obtain
\begin{equation} \label{eq:PPall}
\sum_\pi q^{\vert\pi\vert}=
\prod _{i\ge1} ^{}\frac {1} {(1-q^i)^i} ,
\end{equation}
an elegant product formula for the generating function for
{\it all\/} plane partitions. 

MacMahon developed two very interesting methods to prove
\eqref{eq:PPall} and Theorem~\ref{thm:MM}. 
First, 
MacMahon developed a whole theory, which he called {\it``partition
analysis"}, and which today runs under the name
of {\it``omega calculus}." However, he realised finally that it would
not do what it was supposed to do (namely prove \eqref{eq:PPall}).
So he abandoned this approach. It is interesting to note though that
recently there has been a revival of MacMahon's partition analysis (see
\cite{BrZaAA}), and Andrews and Paule \cite{AnPaAZ}
finally ``made MacMahon's dream true."

The second method consists in translating the problem from 
enumerating plane partitions to the enumeration of --- what
MacMahon called --- {\it lattice permutations}. This idea led
him finally to a proof of Theorem~\ref{thm:MM} in
\cite[Sec.~494]{MacMAA}. The ``correct" general framework for this second
method was brought to light by Richard Stanley in his thesis
(published in revised form as \cite{StanAB}), by developing his
theory of {\it poset partitions}. See the next section, and the article by Ira
Gessel in this volume for an extensive account.

\section{The revival of plane partitions and Richard Stanley}
\label{sec:Stan}

As mentioned in the previous section, 
MacMahon left behind a very intriguing conjecture on ``symmetric"
plane partitions (these are plane partitions which are invariant under
reflection in the main diagonal; see Section~\ref{sec:10}, Class~2).
However, it seems that, at the time, there were not many
others who shared MacMahon's excitement for plane partitions.
(In particular, nobody seemed to care about his conjecture ---
or perhaps it was too difficult at the time \dots)
In any case, after MacMahon plane partitions were more or less forgotten,
except that Wright \cite{WrigAC} calculated the asymptotics of
the number of plane partitions of $n$ as $n$ tends to infinity.
It was Leonard Carlitz \cite{CarlAG,CarlAH}, and 
Basil Gordon and Lorne Houten \cite{GoHoAA,GoHoAB} who,
in the 1960s, relaunched the interest in plane partitions.

However, the ``rebirth" of intensive investigation of plane
partitions was instigated by Stanley's two-part survey article
\cite{StanAA} {\it``Theory and applications of plane partitions"},
together with Bender and Knuth's article \cite{BeKnAA}. 
Already in his thesis (published in revised form as \cite{StanAB}) 
Stanley had dealt with plane
partitions. More precisely, he introduced a vast generalisation
of the notion of plane partitions, {\it poset partitions}, and built
a sophisticated theory around it. In particular, this theory generalised
the earlier mentioned key
idea of MacMahon in his proof of Theorem~\ref{thm:MM} to a 
correspondence between poset partitions and linear extensions
of the underlying poset (again, see the article by Ira Gessel).
In \cite{StanAA}, he laid out the state of affairs in the theory
of plane partitions as of 1971, and he
presented {\it his} view of the subject matter. This meant, among others,
to emphasise the intimate connection with the theory of symmetric
functions. Most importantly, the article \cite{StanAA} offered a
truly fascinating reading, in particular pointing out that
there were several intriguing open problems, some of them 
waiting for a solution already for a very long time.
I am saying all this frankly admitting
  that it was this article of Stanley which ``sparked the fire" 
for plane partitions in me (which is still ``burning" in the form
of my fascination for the enumeration of rhombus tilings).

One of the original contributions one finds in \cite{StanAA} is
an extremely elegant bijective proof of MacMahon's formula
\eqref{eq:PPall}, reported from \cite{StanAC} (but which was published later
than \cite{StanAA}). It is an adaption of an idea of Bender and Knuth
\cite[proof of Theorem~2]{BeKnAA}.
It had the advantage of allowing for an
additional parameter to be put in, the {\it trace} of a plane
partition, which by definition is the sum of the parts along the
diagonal of the plane partition. More on this is to be found in
Section~\ref{sec:trace}.

With the increased interest in the enumeration of plane partition,
several new authors entered the subject, the most prominent being
George Andrews and Ian Macdonald. This led on the one hand to
a proof of MacMahon's conjecture by Andrews \cite{AndrAK},
a proof of a conjecture of Bender and Knuth in the aforementioned
article \cite[p.~50]{BeKnAA} by Gordon \cite{GordAC} and later
by Andrews \cite{AndrAJ}, and alternative proofs of both by
Macdonald \cite[Ex.~16, 17, 19, pp.~83--86]{MacdAC}.
Moreover, Macdonald introduced another symmetry operation on plane
partitions, {\it cyclic symmetry} (see Section~\ref{sec:10}, Class~3),
and observed that, again, the generating function for plane partitions
contained in a given box which are invariant under this operation
seems to be given by an elegant closed form product formula
(see \eqref{eq:C3}).

\begin{figure}[h]
$$
\begin{pmatrix} 
0&0&\hphantom{-}1&\hphantom{-ß}0&0&0\\
1&0&-1&\hphantom{-}1&0&0\\
0&0&\hphantom{-}1&-1&0&1\\
0&1&-1&\hphantom{-}1&0&0\\
0&0&\hphantom{-}1&-1&1&0\\
0&0&\hphantom{-}0&\hphantom{-}1&0&0
\end{pmatrix}
$$
\caption{An alternating sign matrix}
\label{fig:ASM}
\end{figure}

Then, William Mills, David Robbins and Howard Rumsey ``entered the
scene," and Stanley had again a part in this development. Robbins and
Rumsey \cite{RoRuAA} had played with a parametric generalisation of
the determinant which they called {\it$\la$-determinant},\footnote{It
  is interesting to note that this --- at the time, isolated,
and maybe somewhat obscure ---
object has now become part of a fascinating theory, namely 
the theory of discrete integral systems
(cf.\ e.g.\ \cite{DiFrAZ}), and is
maybe the first (non-trivial) example of the so-called {\it Laurent phenomenon}
of Fomin and Zelevinsky \cite{FSZeAF}.} 
and in that context were led to new combinatorial
objects which they called {\it alternating sign matrices}.
An alternating sign matrix is a square matrix consisting of
$0$'s, $1$'s and $(-1)$'s such that, ignoring $0$'s, 
along each row and each column one reads $1,-1,1,\dots,-1,1$
(that is, $1$'s and $(-1)$'s alternate, and at the beginning and at the
end there stands a~$1$). Figure~\ref{fig:ASM} shows a $6\times 6$
alternating sign matrix.
At some point, they became interested in how many such matrices
there were. More precisely, they asked the question: 
how many $n\times n$ alternating sign
matrices are there? 
By computer experiments combined with human
insight, they observed that, apparently, that number was given by
the amazingly simple product formula in \eqref{eq:C10}. 
As Robbins writes in \cite{RobbAA}, after having in vain tried for
some time to prove their conjecture, they began to suspect that this
may have something to do with the theory of plane partitions.
So they called {\it the} expert on plane partitions, Richard Stanley.
After a few days, Stanley replied that, while he was not able to prove
the conjecture, he had seen these numbers before in \cite{AndrAO}, 
namely as the
numbers of ``size~$n$" {\it descending plane partitions}, a certain
variation of cyclically symmetric plane partitions, which Andrews had
introduced in his attempts to prove Macdonald's conjectured formula
for the latter plane partitions. (Andrews ``only" succeeded to prove
the $q=1$ case of Macdonald's conjecture.) 

This led to further remarkable developments in the theory of
plane partitions. 
Andrews and Robbins independently came up with
a conjectured formula for a certain generating function for plane partitions
which are invariant under both reflective and cyclic symmetry (see
Section~\ref{sec:10}, Class~4). 
Moreover, Mills, Robbins and Rumsey now started to look at
alternating sign matrices and plane partitions in parallel, trying
to see connections between these seemingly so different objects.
One of the directions they were investigating was the study of operations on
these objects and whether there existed ``analogous" operations
on the ``other side." One of the outcomes was the discovery
of a new (but, in retrospective, ``obvious") symmetry operation 
on plane partitions: {\it complementation} (see the more detailed
account of this discovery in \cite{RobbAA}). If combined with the
other two symmetry operations --- reflective symmetry, cyclic symmetry
---, this gave rise to six further symmetry classes, beyond the already
``existing" four (see Section~\ref{sec:10}, Classes~5--10). Amazingly,
also for the six new classes, it seemed that the number of plane
partitions in each of these six symmetry classes, the plane partitions
being contained in a given box, was again given by a nice closed form
product formula. Mills, Robbins and Rumsey found striking,
much finer enumerative coincidences between plane partitions
and alternating sign matrices (see Section~\ref{sec:endofstory}), 
and they initiated a
programme of enumerating symmetry classes of alternating sign matrices
(in analogy with the programme for plane partitions summarised here in
Section~\ref{sec:10}). Ironically, they never succeeded to prove any of their
conjectures on alternating sign matrices, but their investigations did
enable them to provide a full proof \cite{MiRRAA} of Macdonald's conjecture on
cyclically symmetric plane partitions and a proof
of the conjectured formula for another symmetry class of plane
partitions (see Section~\ref{sec:10}, Class~8).

In 1986, Stanley summarised the statements, conjectures, and results
(at the time)
for symmetry classes of plane partitions in \cite{StanAI}.
In another survey article \cite{StanAH} he also included the statements,
conjectures, and results for symmetry classes of alternating sign
matrices. So-to-speak, he ``used the opportunity" in \cite{StanAI} 
to resolve one of the (then) conjectures for symmetry
classes of plane partitions, namely the conjecture on the number
of self-complementary plane partitions contained in a given box
(see Section~\ref{sec:10}, Class~6).
Stanley's insight that led him to the solution was that the problem
could be formulated as the problem of calculating a certain
complicated sum, which he identified as a sum of specialised
Schur functions. As Stanley demonstrates, ``by coincidence," this
sum of Schur functions is precisely the expansion of a 
product of two special Schur functions, as an application
of the Littlewood--Richardson rule shows; see Section~\ref{sec:SCPP}.

So, in brief, here is a summary of Stanley's contributions to the
theory of plane partitions:

\begin{itemize} 
\item {\it poset partitions} \cite{StanAB};
\item {\it trace generating functions} for plane partitions
  \cite{StanAC};
\item survey of the theory of plane partitions as of 1971, \cite{StanAA};
\item application of {\it symmetric functions} to the theory of plane
  partitions \cite{StanAA,StanAI};
\item surveys of the state of the art concerning {\it symmetry
classes} of plane partitions and alternating sign matrices as of 
1986 \cite{StanAI,StanAH};
\item proof of the conjecture on {\it self-complementary plane
  partitions} \cite{StanAI}.
\end{itemize}

However, Stanley's contributions extend beyond this through the
work of his students, particularly the work of Emden Gansner
(who continued Stanley's work on traces of plane partitions), Ira
Gessel (who, together with Viennot, developed the powerful method of
non-intersecting lattice paths),
Bob Proctor (who proved the enumeration formulae for Classes~6 and 7
and provided new proofs of several others), 
and John Stembridge (who, among other things, 
paved the way for Andrews' proof of the
enumeration formula for Class~10, and who proved the $q=1$ case of the
enumeration formula for Class~4). Their work will also be
discussed below.

In the next section, I give a quick survey of methods that have
been used to enumerate plane partitions, covering in particular
the non-intersecting path method. Then, in
Section~\ref{sec:trace}, I shall describe Stanley's beautiful bijective
proof of the formula \eqref{eq:PPall} for the generating function for
all plane partitions, leading to a refined formula which features the
statistic ``trace." This section also contains a brief exposition of 
the subsequent work of Gansner on trace generating functions.
Section~\ref{sec:10} presents the project of
enumeration of symmetry classes of plane partitions, listing all
cases and those who proved the corresponding formulae. Pointers to
further work are also included. This is followed by a detailed
discussion of one of the classes in Section~\ref{sec:SCPP}, namely the class of
self-complementary plane partitions. The central piece of that section
is Stanley's proof of the corresponding enumeration formulae.
The final section addresses the role of plane partitions in
research work of today.

\section{Methods for the enumeration of plane partitions}
\label{sec:methods}

It was said in Section~\ref{sec:MM} that plane partitions can be
approached from many different angles. Which are the methods which
have (so far) been applied for the enumeration of plane partitions?

\medskip
As already explained, MacMahon himself already introduced two
methods: first, his {\it partition analysis} (nowadays often called
{\it``omega calculus"}; cf.\ \cite{BrZaAA}), 
which is a generating function method
that, in its utmost generality, can be applied for the enumeration of
integer points in $d$-dimensional space obeying linear inequalities
and equalities.

\medskip
Second, MacMahon introduced a translation of counting problems for
plane partitions to counting problems of so-called ``lattice
permutations," which Stanley \cite{StanAB} generalised to his theory of
{\it poset partitions} (see again Gessel's article in the same volume).

\medskip
{\it The} standard method nowadays for the enumeration of plane
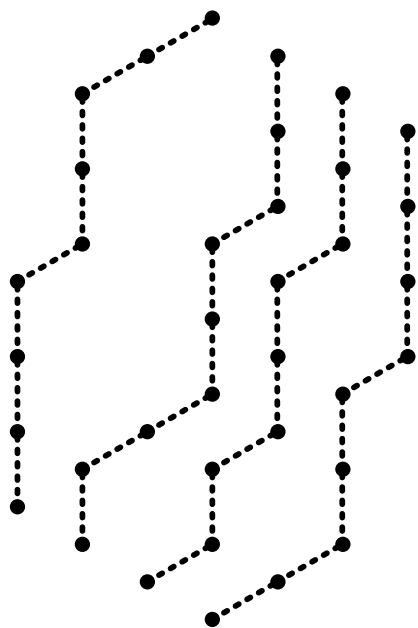
partitions is to use {\it non-intersecting lattice paths} to ``reduce" the
counting problem to the problem of evaluating a certain
determinant.\footnote{Here, ``reduce" is in quotes since, in the harder
cases, the most difficult part then is the evaluation of the
determinant.} I want to illustrate this by explaining the
non-intersecting paths approach to MacMahon's formula in
Theorem~\ref{thm:MM}. For the sake of simplicity, I concentrate on the
$q=1$ case, that is, on the plain enumeration of plane partitions
contained in a given box. The case of generic~$q$ is by no means
more complicated or difficult, it would only require more notation,
which I want to avoid here.

\begin{figure}[h]
\centertexdraw{
  \drawdim truecm  \linewd.02
\move(0 0)
\bsegment
  \drawdim truecm  \linewd.02
  \linewd.12
  \move(0 0)
  \RhombusA \RhombusB \RhombusB 
  \RhombusA \RhombusA \RhombusB \RhombusA \RhombusB \RhombusB
  \move (-0.866025403784439 -.5)
  \RhombusA \RhombusB \RhombusB \RhombusB \RhombusB
  \RhombusA \RhombusA \RhombusB \RhombusA 
  \move (-1.732050807568877 -1)
  \RhombusB \RhombusB \RhombusA \RhombusB \RhombusB \RhombusA
  \RhombusB \RhombusA \RhombusA 
  \move (1.732050807568877 0)
  \RhombusC \RhombusC \RhombusC 
  \move (1.732050807568877 -1)
  \RhombusC \RhombusC \RhombusC 
  \move (3.464101615137755 -3)
  \RhombusC 
  \move (-0.866025403784439 -.5)
  \RhombusC
  \move (-0.866025403784439 -1.5)
  \RhombusC
  \move (0.866025403784439 -2.5)
  \RhombusC \RhombusC 
  \move (0.866025403784439 -3.5)
  \RhombusC \RhombusC \RhombusC 
  \move (2.598076211353316 -5.5)
  \RhombusC 
  \move (0.866025403784439 -5.5)
  \RhombusC 
  \move (-1.732050807568877 -3)
  \RhombusC 
  \move (-1.732050807568877 -4)
  \RhombusC 
  \move (-1.732050807568877 -5)
  \RhombusC \RhombusC 
\linewd.08
\move(-1.29904 -6.25)
\vdSchritt \vdSchritt \vdSchritt 
\hdSchritt \vdSchritt \vdSchritt \hdSchritt \hdSchritt 
\move(-.433012 -6.75)
\vdSchritt \hdSchritt \hdSchritt 
\vdSchritt \vdSchritt \hdSchritt \vdSchritt \vdSchritt 
\move(.433012 -7.25)
\hdSchritt \vdSchritt \hdSchritt 
\vdSchritt \vdSchritt \hdSchritt \vdSchritt \vdSchritt 
\move(1.29904 -7.75)
\hdSchritt \hdSchritt \vdSchritt 
\vdSchritt \hdSchritt \vdSchritt \vdSchritt \vdSchritt 
\ringerl(-1.29904 -6.25)
\ringerl(-.433012 -6.75)
\ringerl(.433012 -7.25)
\ringerl(1.29904 -7.75)
\ringerl(1.29904 .25)
\ringerl(2.165065 -.25)
\ringerl(3.031090 -.75)
\ringerl(3.897115 -1.25)
\esegment
\move(8 0)
\bsegment
  \drawdim truecm  \linewd.02
\linewd.08
\move(-1.29904 -6.25)
\vdSchritt \vdSchritt \vdSchritt 
\hdSchritt \vdSchritt \vdSchritt \hdSchritt \hdSchritt 
\move(-.433012 -6.75)
\vdSchritt \hdSchritt \hdSchritt 
\vdSchritt \vdSchritt \hdSchritt \vdSchritt \vdSchritt 
\move(.433012 -7.25)
\hdSchritt \vdSchritt \hdSchritt 
\vdSchritt \vdSchritt \hdSchritt \vdSchritt \vdSchritt 
\move(1.29904 -7.75)
\hdSchritt \hdSchritt \vdSchritt 
\vdSchritt \hdSchritt \vdSchritt \vdSchritt \vdSchritt 
\ringerl(-1.29904 -6.25)
\ringerl(-1.29904 -5.25)
\ringerl(-1.29904 -4.25)
\ringerl(-1.29904 -3.25)
\ringerl(-.433012 -2.75)
\ringerl(-.433012 -1.75)
\ringerl(-.433012 -.75)
\ringerl(.433012 -.25)
\ringerl(-.433012 -6.75)
\ringerl(-.433012 -5.75)
\ringerl(.433012 -5.25)
\ringerl(1.29904 -4.75)
\ringerl(1.29904 -3.75)
\ringerl(1.29904 -2.75)
\ringerl(2.165065 -2.25)
\ringerl(2.165065 -1.25)
\ringerl(.433012 -7.25)
\ringerl(1.29904 -6.75)
\ringerl(1.29904 -5.75)
\ringerl(2.165065 -5.25)
\ringerl(2.165065 -4.25)
\ringerl(2.165065 -3.25)
\ringerl(3.031090 -2.75)
\ringerl(3.031090 -1.75)
\ringerl(1.29904 -7.75)
\ringerl(2.165065 -7.25)
\ringerl(3.031090 -6.75)
\ringerl(3.031090 -5.75)
\ringerl(3.031090 -4.75)
\ringerl(3.897115 -4.25)
\ringerl(3.897115 -3.25)
\ringerl(3.897115 -2.25)
\ringerl(1.29904 .25)
\ringerl(2.165065 -.25)
\ringerl(3.031090 -.75)
\ringerl(3.897115 -1.25)
\esegment
}
\vskip10pt
\centerline{a. A plane partition\hskip4cm
b. Non-intersecting lattice paths}
$$
\Gitter(4,11)(-5,-0)
\Koordinatenachsen(4,11)(-5,-0)
\SPfad(-4,4),22212211\endSPfad
\SPfad(-3,3),21122122\endSPfad
\SPfad(-2,2),12122122\endSPfad
\SPfad(-1,1),11221222\endSPfad
\DickPunkt(-4,4)
\DickPunkt(-3,3)
\DickPunkt(-2,2)
\DickPunkt(-1,1)
\DickPunkt(2,6)
\DickPunkt(1,7)
\DickPunkt(0,8)
\DickPunkt(-1,9)
$$
\centerline{c. Non-intersecting lattice paths made orthogonal}
\caption{}
\label{fig:PP8}
\end{figure}

Consider the plane partition in
Figure~\ref{fig:PP3}, which is copied in Figure~\ref{fig:PP8}.a. 
We mark the
midpoints of the edges along the south-west side of the hexagon
and we start paths there, where the individual steps of the paths
always connect midpoints of opposite sides of rhombi. See
Figure~\ref{fig:PP8}.a,b
for the result in our running example. We obtain a collection of paths
which connect the midpoints of the edges along the south-west side
with the midpoints of edges along the north-east side. Clearly, the
paths are {\it non-intersecting}, meaning that no two paths have any
points in common. 
By slightly deforming the obtained paths, we may place them into the
plane integer lattice, see Figure~\ref{fig:PP8}.c.
What we have obtained is a bijection between plane partitions in an
$a\times b\times c$ box and families $(P_1,P_2,\dots,P_b)$ of
non-intersecting lattice paths consisting of unit horizontal and
vertical steps in the positive direction, where $P_i$ connects
$(-i,i)$ with $(a-i,c+i)$, $i=1,2,\dots,b$.

The (first) main theorem on non-intersecting lattice paths implies
that the number of the families of non-intersecting lattice paths that
we have obtained above from our plane partitions can be written in
terms of a determinant, see \eqref{eq:PPdet} below. 

\begin{theorem} \label{thm:Lind}%
Let $G$ by a directed, acyclic graph, and
let $A_1,A_2,\dots,A_n$ and\/ $E_1,E_2,\dots,E_n$ be
vertices in $G$ with the property that any path from $A_i$ to $E_l$
and any path from $A_j$ to $E_k$ with $i<j$ and $k<l$
have a common vertex. Then the number of all 
families $(P_1,P_2,\dots,P_n)$ of non-intersecting 
paths, where $P_i$ runs from $A_i$ to $E_i$, $i=1,2,\dots,n$, 
is given by
\begin{equation} \label{CKe2.5a}
\det_{1\le i,j\le n}\big(L_G(A_j\to E_i)\big),
\end{equation}
where $L_G(A\to E)$ denotes the number of all paths 
starting in $A$ and ending in $E$.
\end{theorem}

This theorem was discovered and independently rediscovered several times. 
It is originally due to Lindstr\"om \cite[Lemma~1]{LindAA}
(who, in fact, proved a more general theorem, and which also
includes weights in a straightforward way; for a proof of
Theorem~\ref{thm:MM} with generic~$q$, we would need the weighted
version), who discovered and used it in the context of matroid
representations. It was Gessel and Viennot \cite{GeViAA,GeViAB} who realised
its enormous significance for the enumeration of plane 
partitions.\footnote{Lindstr\"om's theorem was rediscovered
(not always in its most general form) in the 1980s at about the same
time in three different
communities, not knowing of each other at that time: 
in enumerative combinatorics
by Gessel and Viennot \cite{GeViAA,GeViAB} in order to count
tableaux and plane partitions, in statistical
physics by Fisher \cite[Sec.~5.3]{FishAA} in order to apply it to
the analysis of vicious walkers as a model of wetting and melting, 
and in combinatorial chemistry by 
John and Sachs \cite{JoSaAB} and Gronau, Just, Schade, Scheffler and 
Wojciechowski \cite{GrJSAA} in order to compute Pauling's bond order
in benzenoid hydrocarbon molecules. 
It must however be
mentioned that in fact the same idea appeared even earlier in work by
Karlin and McGregor \cite{KaMGAB,KaMGAC} in a probabilistic
framework, as well as that the so-called 
``Slater determinant"
in quantum mechanics (cf.\ \cite{SlatZY} and \cite[Ch.~11]{SlatZZ}) 
may qualify as an ``ancestor" of the determinantal formula of
Lindstr\"om.}

If we apply Theorem~\ref{thm:Lind} with the graph consisting
of the integer points in the plane as vertices and edges
being horizontal and vertical unit vectors in the positive direction,
with $A_i=(-i,i)$ and $E_i=(a-i,c+i)$, then we obtain that the number
of plane partitions contained in an $a\times b\times c$ box equals the
determinant 
\begin{equation} \label{eq:PPdet}
\det_{1\le i,j\le b}\left(\binom {a+c}{a-i+j}\right). 
\end{equation}
There are many ways to evaluate this determinant, see e.g.\
\cite[Secs.~2.2, 2.3, 2.5]{KratBN}, and the result can be written in the form 
\begin{equation} \label{eq:MM} 
\prod _{i=1} ^{a}
\prod _{j=1} ^{b}
\prod _{k=1} ^{c}
\frac {i+j+k-1} {i+j+k-2},
\end{equation}
which is the $q\to1$ limit of \eqref{eq:M2}.

\medskip
Theorem~\ref{thm:Lind} solves the problem of counting
non-intersecting lattice paths for which starting {\it and\/} end
points are fixed.
If one tries the non-intersecting lattice path approach for
symmetry classes of plane partitions, then one often encounters the
problem that starting or/and end points are not fixed, but rather may
vary in given sets. Inspired by work of Okada \cite{OkadAA} on summations
of minors of given matrices (which was later generalised to
the minor summation theorem of Ishikawa and Wakayama \cite{IsWaAA}), 
Stembridge developed the corresponding theory, which provides
Pfaffian formulae (and thus again determinantal formulae, 
given the fact that the Pfaffian is the square
root of a determinant) 
for the number of non-intersecting lattice paths where either starting
or end points, or both, vary in given sets. With the exception of
only very few, it is the non-intersecting lattice path method
which constitutes the first step in proofs of enumeration formulae
for symmetry classes of plane partitions, may it be explicitly or
implicitly. 
We refer the reader to Section~\ref{sec:SCPP} for 
an ``implicit" application of non-intersecting lattice paths
(in the sense that Stanley's original proof did not mention
non-intersecting lattice paths; they are still there).

\medskip
Symmetric functions and representation theory come into play because
plane partitions are very close to the tableau-like objects that index
representations of classical groups. As said earlier, Stanley
\cite{StanAA} was the
first to emphasise and exploit this connection. It was picked up and
deepened by Macdonald \cite[Ch.~I, Sec.~5, Ex.~13--19]{MacdAC}, Proctor
\cite{ProcAE,ProcAD,ProcAB}, Stembridge
\cite{StemAD,StemAE,StemAI}, Okounkov, Reshetikhin,
and Vuleti\'c (see \cite{VuleAA} and the references therein), 
and Kuperberg \cite{KupeAF}. 
To give a flavour, let us
again consider the enumeration of all plane partitions in
an $a\times b\times c$ box. In terms of the original definition
\eqref{eq:PP}, these are plane partitions of shape
$(b,b,\dots,b)$ (with $a$ occurrences of $b$) with entries at
most~$c$. An example with $a=3$, $b=4$, $c=6$ is shown in
Figure~\ref{fig:PP9}.a. (The corresponding ``graphical"
representations are the ones in Figures~\ref{fig:PP2} and
\ref{fig:PP3}.) We rotate the plane partition by $180^\circ$
and add $i$ to each entry in row $i$, for $i=1,2,\dots,c$.
Thereby we obtain an array of integers of the same shape as the
original plane partition, however with the property that entries along
rows are weakly increasing and entries along columns are strictly
increasing, all of them between $1$ and $a+c$. 
Arrays satisfying the above two monotonicity properties are called {\it
  semistandard tableaux}. See Figure~\ref{fig:PP9}.b for the
semistandard tableau which arises in this way from the plane partition
in Figure~\ref{fig:PP9}.a. 

\begin{figure}[h]
$$
\begin{matrix}
5&3&3&2\\
5&1&1&0\\
3&1&0&0
\end{matrix}
\hskip6cm
\begin{matrix}
1&1&2&4\\
2&3&3&7\\
5&6&6&8
\end{matrix}
\hskip2cm
$$
\vskip10pt
\centerline{a. A plane partition\hskip2cm
b. The corresponding semistandard tableau}
\label{fig:PP9}
\caption{}
\end{figure}

It is a central fact of the representation theory of $SL_n(\mathbb C)$
that semistandard tableaux of shape $\la$ 
(where shape is defined in the same way as for a
plane partition) with entries between $1$ and $n$
index a basis of an irreducible representation of $SL_n(\mathbb C)$.
On the character level, this is reflected by the {\it Schur function}
$s_\la(x_1,x_2,\dots,x_n)$ (cf.\ \cite{MacdAC}), which is defined by
\begin{equation} \label{eq:Schur} 
s_\la(x_1,x_2,\dots,x_n)=\sum_T 
\prod _{i=1} ^{n}x_i^{\#\text{entries $i$ in T}},
\end{equation}
where the sum is over all semistandard tableaux $T$ of shape $\la$.
Thus, we see that MacMahon's generating function for plane partitions
essentially equals a specialised Schur function, namely we have
\begin{equation} \label{eq:MMSchur}
\sum_\pi q^{\vert\pi\vert}=q^{-b\binom {a+1}2}
s_{(b,b,\dots,b)}(q,q^2,\dots,q^{a+c}),
\end{equation}
the sum on the left-hand side being over all plane partitions
$\pi$ contained in an $a\times b\times c$ box.
If this ``principal specialisation" of the variables,
i.e., $x_i=q^i$, $i=1,2,\dots,a+c$, is
done in the Weyl character formula for Schur functions
(cf.\ \cite[p.~40, Eq.~(3.1)]{MacdAC}),
$$
s_\la(x_1,x_2,\dots,x_n)=\frac {\det_{1\le i,j\le n}(x_i^{\la_j+n-j})} 
{\det_{1\le i,j\le n}(x_i^{n-j})} ,
$$
then both determinants are Vandermonde determinants and can therefore
be evaluated, thus establishing \eqref{eq:M2}.

\begin{figure}[h]
\centertexdraw{
  \drawdim truecm  \linewd.02
  \rhombus \rhombus \rhombus \rhombus \ldreieck
  \move (-0.866025403784439 -.5)
  \rhombus \rhombus \rhombus \rhombus \rhombus \ldreieck
  \move (-1.732050807568877 -1)
  \rhombus \rhombus \rhombus \rhombus \rhombus \rhombus \ldreieck
  \move (-1.732050807568877 -1)
  \rdreieck
  \rhombus \rhombus \rhombus \rhombus \rhombus \rhombus \ldreieck
  \move (-1.732050807568877 -2)
  \rdreieck
  \rhombus \rhombus \rhombus \rhombus \rhombus \rhombus \ldreieck
  \move (-1.732050807568877 -3)
  \rdreieck
  \rhombus \rhombus \rhombus \rhombus \rhombus \rhombus 
  \move (-1.732050807568877 -4)
  \rdreieck
  \rhombus \rhombus \rhombus \rhombus \rhombus 
  \move (-1.732050807568877 -5)
  \rdreieck
  \rhombus \rhombus \rhombus \rhombus 
  \linewd.05
  \move (0.57735 0)
  \hhSchritt \rdSchritt 
  \hhSchritt \rdSchritt 
  \hhSchritt \rdSchritt 
  \hhSchritt \rdSchritt 
  \ldSchritt \rdSchritt 
  \ldSchritt \rdSchritt 
  \ldSchritt \rdSchritt 
  \ldSchritt \rdSchritt 
  \ldSchritt \hhhSchritt 
  \ldSchritt \hhhSchritt 
  \ldSchritt \hhhSchritt 
  \move (0.57735 0)
  \ldSchritt 
  \hhhSchritt \ldSchritt 
  \hhhSchritt \ldSchritt 
  \rdSchritt \ldSchritt 
  \rdSchritt \ldSchritt 
  \rdSchritt \ldSchritt 
  \rdSchritt \ldSchritt 
  \rdSchritt \hhSchritt 
  \rdSchritt \hhSchritt 
  \rdSchritt \hhSchritt 
  \rdSchritt \hhSchritt 
  \move (0.57735 -1)
  \hhSchritt \rdSchritt 
  \hhSchritt \rdSchritt 
  \hhSchritt \rdSchritt 
  \ldSchritt \rdSchritt 
  \ldSchritt \rdSchritt 
  \ldSchritt \rdSchritt 
  \ldSchritt \hhhSchritt 
  \ldSchritt \hhhSchritt 
  \move (0.57735 -1)
  \ldSchritt 
  \hhhSchritt \ldSchritt 
  \rdSchritt \ldSchritt 
  \rdSchritt \ldSchritt 
  \rdSchritt \ldSchritt 
  \rdSchritt \hhSchritt 
  \rdSchritt \hhSchritt 
  \rdSchritt \hhSchritt 
  \move (0.57735 -2)
  \hhSchritt \rdSchritt 
  \hhSchritt \rdSchritt 
  \ldSchritt \rdSchritt 
  \ldSchritt \rdSchritt 
  \ldSchritt \hhhSchritt 
  \move (0.57735 -2)
  \ldSchritt 
  \rdSchritt \ldSchritt 
  \rdSchritt \ldSchritt 
  \rdSchritt \hhSchritt 
  \rdSchritt \hhSchritt 
  \move (0.57735 -3)
  \hhSchritt \rdSchritt 
  \ldSchritt \rdSchritt 
  \ldSchritt 
  \move (0.288675 -.5)
  \rdSchritt 
  \move (-0.57735 -1)
  \rdSchritt 
  \move (1.44338 -.5)
  \ldSchritt 
  \move (2.3094 -1)
  \ldSchritt 
  \move (3.17543 -1.5)
  \ldSchritt 
  \move (0.288675 -1.5)
  \rdSchritt 
  \move (1.44338 -1.5)
  \ldSchritt 
  \move (2.3094 -2)
  \ldSchritt 
  \move (-1.1547 -2)
  \hhSchritt 
  \move (-1.1547 -3)
  \hhSchritt 
  \move (-1.1547 -4)
  \hhSchritt 
  \move (-1.1547 -5)
  \hhSchritt 
  \move (-0.288675 -2.5)
  \hhSchritt 
  \move (-0.288675 -3.5)
  \hhSchritt 
  \move (-0.288675 -4.5)
  \hhSchritt 
  \move (0.57735 -4)
  \hhSchritt 
  \move (1.44338 -3.5)
  \hhSchritt 
  \move (1.44338 -4.5)
  \hhSchritt 
  \move (2.3094 -3)
  \hhSchritt 
  \move (2.3094 -4)
  \hhSchritt 
  \move (2.3094 -5)
  \hhSchritt 
  \move (3.17543 -2.5)
  \hhSchritt 
  \move (3.17543 -3.5)
  \hhSchritt 
  \move (3.17543 -4.5)
  \hhSchritt 
  \move (3.17543 -5.5)
  \hhSchritt 
  \move (1.1547 -3)
  \ruSchritt 
  \move (1.1547 -6)
  \ruSchritt 
  \move (1.1547 -7)
  \ruSchritt 
  \move (0.288675 -5.5)
  \ruSchritt 
  \move (0.288675 -6.5)
  \ruSchritt 
  \move (-0.57735 -6)
  \ruSchritt 
  \move (2.3094 -6)
  \luSchritt 
  \move (2.3094 -7)
  \luSchritt 
  \move (3.17543 -6.5)
  \luSchritt 
  \ringerl (0.57735 0)
  \ringerl (0.57735 -1)
  \ringerl (0.57735 -2)
  \ringerl (0.57735 -3)
  \ringerl (0.57735 -4)
  \ringerl (0.57735 -5)
  \ringerl (0.57735 -6)
  \ringerl (0.57735 -7)
  \ringerl (-0.288675 -.5)
  \ringerl (-0.288675 -1.5)
  \ringerl (-0.288675 -2.5)
  \ringerl (-0.288675 -3.5)
  \ringerl (-0.288675 -4.5)
  \ringerl (-0.288675 -5.5)
  \ringerl (-0.288675 -6.5)
  \ringerl (-1.1547 -1)
  \ringerl (-1.1547 -2)
  \ringerl (-1.1547 -3)
  \ringerl (-1.1547 -4)
  \ringerl (-1.1547 -5)
  \ringerl (-1.1547 -6)
  \ringerl (1.1547 0)
  \ringerl (1.1547 -1)
  \ringerl (1.1547 -2)
  \ringerl (1.1547 -3)
  \ringerl (1.1547 -4)
  \ringerl (1.1547 -5)
  \ringerl (1.1547 -6)
  \ringerl (1.1547 -7)
  \ringerl (2.02073 -.5)
  \ringerl (2.02073 -1.5)
  \ringerl (2.02073 -2.5)
  \ringerl (2.02073 -3.5)
  \ringerl (2.02073 -4.5)
  \ringerl (2.02073 -5.5)
  \ringerl (2.02073 -6.5)
  \ringerl (2.02073 -7.5)
  \ringerl (2.88675 -1)
  \ringerl (2.88675 -2)
  \ringerl (2.88675 -3)
  \ringerl (2.88675 -4)
  \ringerl (2.88675 -5)
  \ringerl (2.88675 -6)
  \ringerl (2.88675 -7)
  \ringerl (3.75278 -1.5)
  \ringerl (3.75278 -2.5)
  \ringerl (3.75278 -3.5)
  \ringerl (3.75278 -4.5)
  \ringerl (3.75278 -5.5)
  \ringerl (3.75278 -6.5)
  \ringerl (-1.44338 -1.5)
  \ringerl (-1.44338 -2.5)
  \ringerl (-1.44338 -3.5)
  \ringerl (-1.44338 -4.5)
  \ringerl (-1.44338 -5.5)
  \ringerl (-0.57735 -1)
  \ringerl (-0.57735 -2)
  \ringerl (-0.57735 -3)
  \ringerl (-0.57735 -4)
  \ringerl (-0.57735 -5)
  \ringerl (-0.57735 -6)
  \ringerl (0.288675 -.5)
  \ringerl (0.288675 -1.5)
  \ringerl (0.288675 -2.5)
  \ringerl (0.288675 -3.5)
  \ringerl (0.288675 -4.5)
  \ringerl (0.288675 -5.5)
  \ringerl (0.288675 -6.5)
  \ringerl (1.44338 -.5)
  \ringerl (1.44338 -1.5)
  \ringerl (1.44338 -2.5)
  \ringerl (1.44338 -3.5)
  \ringerl (1.44338 -4.5)
  \ringerl (1.44338 -5.5)
  \ringerl (1.44338 -6.5)
  \ringerl (1.44338 -7.5)
  \ringerl (2.3094 -1)
  \ringerl (2.3094 -2)
  \ringerl (2.3094 -3)
  \ringerl (2.3094 -4)
  \ringerl (2.3094 -5)
  \ringerl (2.3094 -6)
  \ringerl (2.3094 -7)
  \ringerl (3.17543 -1.5)
  \ringerl (3.17543 -2.5)
  \ringerl (3.17543 -3.5)
  \ringerl (3.17543 -4.5)
  \ringerl (3.17543 -5.5)
  \ringerl (3.17543 -6.5)
  \ringerl (4.04145 -2)
  \ringerl (4.04145 -3)
  \ringerl (4.04145 -4)
  \ringerl (4.04145 -5)
  \ringerl (4.04145 -6)
\move(8 0)
\bsegment
  \drawdim truecm  \linewd.02
  \ringerl (0.57735 0)
  \ringerl (0.57735 -1)
  \ringerl (0.57735 -2)
  \ringerl (0.57735 -3)
  \ringerl (0.57735 -4)
  \ringerl (0.57735 -5)
  \ringerl (0.57735 -6)
  \ringerl (0.57735 -7)
  \ringerl (-0.288675 -.5)
  \ringerl (-0.288675 -1.5)
  \ringerl (-0.288675 -2.5)
  \ringerl (-0.288675 -3.5)
  \ringerl (-0.288675 -4.5)
  \ringerl (-0.288675 -5.5)
  \ringerl (-0.288675 -6.5)
  \ringerl (-1.1547 -1)
  \ringerl (-1.1547 -2)
  \ringerl (-1.1547 -3)
  \ringerl (-1.1547 -4)
  \ringerl (-1.1547 -5)
  \ringerl (-1.1547 -6)
  \ringerl (1.1547 0)
  \ringerl (1.1547 -1)
  \ringerl (1.1547 -2)
  \ringerl (1.1547 -3)
  \ringerl (1.1547 -4)
  \ringerl (1.1547 -5)
  \ringerl (1.1547 -6)
  \ringerl (1.1547 -7)
  \ringerl (2.02073 -.5)
  \ringerl (2.02073 -1.5)
  \ringerl (2.02073 -2.5)
  \ringerl (2.02073 -3.5)
  \ringerl (2.02073 -4.5)
  \ringerl (2.02073 -5.5)
  \ringerl (2.02073 -6.5)
  \ringerl (2.02073 -7.5)
  \ringerl (2.88675 -1)
  \ringerl (2.88675 -2)
  \ringerl (2.88675 -3)
  \ringerl (2.88675 -4)
  \ringerl (2.88675 -5)
  \ringerl (2.88675 -6)
  \ringerl (2.88675 -7)
  \ringerl (3.75278 -1.5)
  \ringerl (3.75278 -2.5)
  \ringerl (3.75278 -3.5)
  \ringerl (3.75278 -4.5)
  \ringerl (3.75278 -5.5)
  \ringerl (3.75278 -6.5)
  \ringerl (-1.44338 -1.5)
  \ringerl (-1.44338 -2.5)
  \ringerl (-1.44338 -3.5)
  \ringerl (-1.44338 -4.5)
  \ringerl (-1.44338 -5.5)
  \ringerl (-0.57735 -1)
  \ringerl (-0.57735 -2)
  \ringerl (-0.57735 -3)
  \ringerl (-0.57735 -4)
  \ringerl (-0.57735 -5)
  \ringerl (-0.57735 -6)
  \ringerl (0.288675 -.5)
  \ringerl (0.288675 -1.5)
  \ringerl (0.288675 -2.5)
  \ringerl (0.288675 -3.5)
  \ringerl (0.288675 -4.5)
  \ringerl (0.288675 -5.5)
  \ringerl (0.288675 -6.5)
  \ringerl (1.44338 -.5)
  \ringerl (1.44338 -1.5)
  \ringerl (1.44338 -2.5)
  \ringerl (1.44338 -3.5)
  \ringerl (1.44338 -4.5)
  \ringerl (1.44338 -5.5)
  \ringerl (1.44338 -6.5)
  \ringerl (1.44338 -7.5)
  \ringerl (2.3094 -1)
  \ringerl (2.3094 -2)
  \ringerl (2.3094 -3)
  \ringerl (2.3094 -4)
  \ringerl (2.3094 -5)
  \ringerl (2.3094 -6)
  \ringerl (2.3094 -7)
  \ringerl (3.17543 -1.5)
  \ringerl (3.17543 -2.5)
  \ringerl (3.17543 -3.5)
  \ringerl (3.17543 -4.5)
  \ringerl (3.17543 -5.5)
  \ringerl (3.17543 -6.5)
  \ringerl (4.04145 -2)
  \ringerl (4.04145 -3)
  \ringerl (4.04145 -4)
  \ringerl (4.04145 -5)
  \ringerl (4.04145 -6)
  \linewd.12
  \move(0.57735 0)
  \HHSchritt
  \move(-0.288675 -.5)
  \HHSchritt
  \move(-1.1547 -1)
  \LDSchritt
  \move(-0.57735 -1)
  \RDSchritt
  \move(0.57735 -1)
  \LDSchritt
  \move(1.1547 -1)
  \RUSchritt
  \move(2.02073 -.5)
  \RDSchritt
  \move(2.88675 -1)
  \RDSchritt
  \move(3.75278 -1.5)
  \RDSchritt
  \move(-1.1547 -2)
  \LDSchritt
  \move(-0.57735 -2)
  \RDSchritt
  \move(0.57735 -2)
  \LDSchritt
  \move(1.1547 -2)
  \RUSchritt
  \move(2.02073 -1.5)
  \RDSchritt
  \move(2.88675 -2)
  \RDSchritt
  \move(3.75278 -2.5)
  \RDSchritt
  \move(-1.1547 -3)
  \HHSchritt
  \move(1.44338 -2.5)
  \HHSchritt
  \move(2.3094 -3)
  \HHSchritt
  \move(-1.1547 -4)
  \LUSchritt
  \move(-0.57735 -4)
  \RUSchritt
  \move(0.57735 -3)
  \LDSchritt
  \move(1.1547 -3)
  \RDSchritt
  \move(2.02073 -3.5)
  \RDSchritt
  \move(2.88675 -4)
  \RUSchritt
  \move(3.75278 -3.5)
  \RDSchritt
  \move(-1.1547 -5)
  \LUSchritt
  \move(-0.57735 -5)
  \RUSchritt
  \move(0.57735 -4)
  \LDSchritt
  \move(1.1547 -4)
  \RDSchritt
  \move(2.02073 -4.5)
  \RDSchritt
  \move(2.88675 -5)
  \RDSchritt
  \move(3.75278 -5.5)
  \RUSchritt
  \move(0.57735 -5)
  \HHSchritt
  \move(-0.288675 -5.5)
  \HHSchritt
  \move(1.44338 -5.5)
  \HHSchritt
  \move(3.17543 -4.5)
  \HHSchritt
  \move(0.57735 -7)
  \HHSchritt
  \move(1.44338 -7.5)
  \HHSchritt
  \move(2.3094 -7)
  \HHSchritt
  \move(-1.1547 -6)
  \LUSchritt
  \move(-0.57735 -6)
  \RDSchritt
  \move(0.57735 -6)
  \LDSchritt
  \move(1.1547 -6)
  \RDSchritt
  \move(2.02073 -6.5)
  \RUSchritt
  \move(2.88675 -6)
  \RDSchritt
  \move(3.75278 -6.5)
  \RUSchritt
  \linewd.08
  \move(0 0)
  \RhombusA \RhombusB \RhombusB 
  \RhombusA \RhombusA \RhombusB \RhombusA \RhombusB \RhombusB
  \move (-0.866025403784439 -.5)
  \RhombusA \RhombusB \RhombusB \RhombusB \RhombusB
  \RhombusA \RhombusA \RhombusB \RhombusA 
  \move (-1.732050807568877 -1)
  \RhombusB \RhombusB \RhombusA \RhombusB \RhombusB \RhombusA
  \RhombusB \RhombusA \RhombusA 
  \move (1.732050807568877 0)
  \RhombusC \RhombusC \RhombusC 
  \move (1.732050807568877 -1)
  \RhombusC \RhombusC \RhombusC 
  \move (3.464101615137755 -3)
  \RhombusC 
  \move (-0.866025403784439 -.5)
  \RhombusC
  \move (-0.866025403784439 -1.5)
  \RhombusC
  \move (0.866025403784439 -2.5)
  \RhombusC \RhombusC 
  \move (0.866025403784439 -3.5)
  \RhombusC \RhombusC \RhombusC 
  \move (2.598076211353316 -5.5)
  \RhombusC 
  \move (0.866025403784439 -5.5)
  \RhombusC 
  \move (-1.732050807568877 -3)
  \RhombusC 
  \move (-1.732050807568877 -4)
  \RhombusC 
  \move (-1.732050807568877 -5)
  \RhombusC \RhombusC 
\esegment
\htext (-2.5 -9){\small a. The triangular grid with its ``dual" graph}
\htext (5.8 -9){\small b. A rhombus tiling and corresponding}
\htext (5.8 -9.5){\small \hphantom{b. }edges in the ``dual" graph}
\move(4 -10)
\bsegment
  \drawdim truecm  \linewd.02
  \linewd.05
  \move (0.57735 0)
  \hhSchritt \rdSchritt 
  \hhSchritt \rdSchritt 
  \hhSchritt \rdSchritt 
  \hhSchritt \rdSchritt 
  \ldSchritt \rdSchritt 
  \ldSchritt \rdSchritt 
  \ldSchritt \rdSchritt 
  \ldSchritt \rdSchritt 
  \ldSchritt \hhhSchritt 
  \ldSchritt \hhhSchritt 
  \ldSchritt \hhhSchritt 
  \move (0.57735 0)
  \ldSchritt 
  \hhhSchritt \ldSchritt 
  \hhhSchritt \ldSchritt 
  \rdSchritt \ldSchritt 
  \rdSchritt \ldSchritt 
  \rdSchritt \ldSchritt 
  \rdSchritt \ldSchritt 
  \rdSchritt \hhSchritt 
  \rdSchritt \hhSchritt 
  \rdSchritt \hhSchritt 
  \rdSchritt \hhSchritt 
  \move (0.57735 -1)
  \hhSchritt \rdSchritt 
  \hhSchritt \rdSchritt 
  \hhSchritt \rdSchritt 
  \ldSchritt \rdSchritt 
  \ldSchritt \rdSchritt 
  \ldSchritt \rdSchritt 
  \ldSchritt \hhhSchritt 
  \ldSchritt \hhhSchritt 
  \move (0.57735 -1)
  \ldSchritt 
  \hhhSchritt \ldSchritt 
  \rdSchritt \ldSchritt 
  \rdSchritt \ldSchritt 
  \rdSchritt \ldSchritt 
  \rdSchritt \hhSchritt 
  \rdSchritt \hhSchritt 
  \rdSchritt \hhSchritt 
  \move (0.57735 -2)
  \hhSchritt \rdSchritt 
  \hhSchritt \rdSchritt 
  \ldSchritt \rdSchritt 
  \ldSchritt \rdSchritt 
  \ldSchritt \hhhSchritt 
  \move (0.57735 -2)
  \ldSchritt 
  \rdSchritt \ldSchritt 
  \rdSchritt \ldSchritt 
  \rdSchritt \hhSchritt 
  \rdSchritt \hhSchritt 
  \move (0.57735 -3)
  \hhSchritt \rdSchritt 
  \ldSchritt \rdSchritt 
  \ldSchritt 
  \move (0.288675 -.5)
  \rdSchritt 
  \move (-0.57735 -1)
  \rdSchritt 
  \move (1.44338 -.5)
  \ldSchritt 
  \move (2.3094 -1)
  \ldSchritt 
  \move (3.17543 -1.5)
  \ldSchritt 
  \move (0.288675 -1.5)
  \rdSchritt 
  \move (1.44338 -1.5)
  \ldSchritt 
  \move (2.3094 -2)
  \ldSchritt 
  \move (-1.1547 -2)
  \hhSchritt 
  \move (-1.1547 -3)
  \hhSchritt 
  \move (-1.1547 -4)
  \hhSchritt 
  \move (-1.1547 -5)
  \hhSchritt 
  \move (-0.288675 -2.5)
  \hhSchritt 
  \move (-0.288675 -3.5)
  \hhSchritt 
  \move (-0.288675 -4.5)
  \hhSchritt 
  \move (0.57735 -4)
  \hhSchritt 
  \move (1.44338 -3.5)
  \hhSchritt 
  \move (1.44338 -4.5)
  \hhSchritt 
  \move (2.3094 -3)
  \hhSchritt 
  \move (2.3094 -4)
  \hhSchritt 
  \move (2.3094 -5)
  \hhSchritt 
  \move (3.17543 -2.5)
  \hhSchritt 
  \move (3.17543 -3.5)
  \hhSchritt 
  \move (3.17543 -4.5)
  \hhSchritt 
  \move (3.17543 -5.5)
  \hhSchritt 
  \move (1.1547 -3)
  \ruSchritt 
  \move (1.1547 -6)
  \ruSchritt 
  \move (1.1547 -7)
  \ruSchritt 
  \move (0.288675 -5.5)
  \ruSchritt 
  \move (0.288675 -6.5)
  \ruSchritt 
  \move (-0.57735 -6)
  \ruSchritt 
  \move (2.3094 -6)
  \luSchritt 
  \move (2.3094 -7)
  \luSchritt 
  \move (3.17543 -6.5)
  \luSchritt 
  \ringerl (0.57735 0)
  \ringerl (0.57735 -1)
  \ringerl (0.57735 -2)
  \ringerl (0.57735 -3)
  \ringerl (0.57735 -4)
  \ringerl (0.57735 -5)
  \ringerl (0.57735 -6)
  \ringerl (0.57735 -7)
  \ringerl (-0.288675 -.5)
  \ringerl (-0.288675 -1.5)
  \ringerl (-0.288675 -2.5)
  \ringerl (-0.288675 -3.5)
  \ringerl (-0.288675 -4.5)
  \ringerl (-0.288675 -5.5)
  \ringerl (-0.288675 -6.5)
  \ringerl (-1.1547 -1)
  \ringerl (-1.1547 -2)
  \ringerl (-1.1547 -3)
  \ringerl (-1.1547 -4)
  \ringerl (-1.1547 -5)
  \ringerl (-1.1547 -6)
  \ringerl (1.1547 0)
  \ringerl (1.1547 -1)
  \ringerl (1.1547 -2)
  \ringerl (1.1547 -3)
  \ringerl (1.1547 -4)
  \ringerl (1.1547 -5)
  \ringerl (1.1547 -6)
  \ringerl (1.1547 -7)
  \ringerl (2.02073 -.5)
  \ringerl (2.02073 -1.5)
  \ringerl (2.02073 -2.5)
  \ringerl (2.02073 -3.5)
  \ringerl (2.02073 -4.5)
  \ringerl (2.02073 -5.5)
  \ringerl (2.02073 -6.5)
  \ringerl (2.02073 -7.5)
  \ringerl (2.88675 -1)
  \ringerl (2.88675 -2)
  \ringerl (2.88675 -3)
  \ringerl (2.88675 -4)
  \ringerl (2.88675 -5)
  \ringerl (2.88675 -6)
  \ringerl (2.88675 -7)
  \ringerl (3.75278 -1.5)
  \ringerl (3.75278 -2.5)
  \ringerl (3.75278 -3.5)
  \ringerl (3.75278 -4.5)
  \ringerl (3.75278 -5.5)
  \ringerl (3.75278 -6.5)
  \ringerl (-1.44338 -1.5)
  \ringerl (-1.44338 -2.5)
  \ringerl (-1.44338 -3.5)
  \ringerl (-1.44338 -4.5)
  \ringerl (-1.44338 -5.5)
  \ringerl (-0.57735 -1)
  \ringerl (-0.57735 -2)
  \ringerl (-0.57735 -3)
  \ringerl (-0.57735 -4)
  \ringerl (-0.57735 -5)
  \ringerl (-0.57735 -6)
  \ringerl (0.288675 -.5)
  \ringerl (0.288675 -1.5)
  \ringerl (0.288675 -2.5)
  \ringerl (0.288675 -3.5)
  \ringerl (0.288675 -4.5)
  \ringerl (0.288675 -5.5)
  \ringerl (0.288675 -6.5)
  \ringerl (1.44338 -.5)
  \ringerl (1.44338 -1.5)
  \ringerl (1.44338 -2.5)
  \ringerl (1.44338 -3.5)
  \ringerl (1.44338 -4.5)
  \ringerl (1.44338 -5.5)
  \ringerl (1.44338 -6.5)
  \ringerl (1.44338 -7.5)
  \ringerl (2.3094 -1)
  \ringerl (2.3094 -2)
  \ringerl (2.3094 -3)
  \ringerl (2.3094 -4)
  \ringerl (2.3094 -5)
  \ringerl (2.3094 -6)
  \ringerl (2.3094 -7)
  \ringerl (3.17543 -1.5)
  \ringerl (3.17543 -2.5)
  \ringerl (3.17543 -3.5)
  \ringerl (3.17543 -4.5)
  \ringerl (3.17543 -5.5)
  \ringerl (3.17543 -6.5)
  \ringerl (4.04145 -2)
  \ringerl (4.04145 -3)
  \ringerl (4.04145 -4)
  \ringerl (4.04145 -5)
  \ringerl (4.04145 -6)
  \linewd.12
  \move(0.57735 0)
  \HHSchritt
  \move(-0.288675 -.5)
  \HHSchritt
  \move(-1.1547 -1)
  \LDSchritt
  \move(-0.57735 -1)
  \RDSchritt
  \move(0.57735 -1)
  \LDSchritt
  \move(1.1547 -1)
  \RUSchritt
  \move(2.02073 -.5)
  \RDSchritt
  \move(2.88675 -1)
  \RDSchritt
  \move(3.75278 -1.5)
  \RDSchritt
  \move(-1.1547 -2)
  \LDSchritt
  \move(-0.57735 -2)
  \RDSchritt
  \move(0.57735 -2)
  \LDSchritt
  \move(1.1547 -2)
  \RUSchritt
  \move(2.02073 -1.5)
  \RDSchritt
  \move(2.88675 -2)
  \RDSchritt
  \move(3.75278 -2.5)
  \RDSchritt
  \move(-1.1547 -3)
  \HHSchritt
  \move(1.44338 -2.5)
  \HHSchritt
  \move(2.3094 -3)
  \HHSchritt
  \move(-1.1547 -4)
  \LUSchritt
  \move(-0.57735 -4)
  \RUSchritt
  \move(0.57735 -3)
  \LDSchritt
  \move(1.1547 -3)
  \RDSchritt
  \move(2.02073 -3.5)
  \RDSchritt
  \move(2.88675 -4)
  \RUSchritt
  \move(3.75278 -3.5)
  \RDSchritt
  \move(-1.1547 -5)
  \LUSchritt
  \move(-0.57735 -5)
  \RUSchritt
  \move(0.57735 -4)
  \LDSchritt
  \move(1.1547 -4)
  \RDSchritt
  \move(2.02073 -4.5)
  \RDSchritt
  \move(2.88675 -5)
  \RDSchritt
  \move(3.75278 -5.5)
  \RUSchritt
  \move(0.57735 -5)
  \HHSchritt
  \move(-0.288675 -5.5)
  \HHSchritt
  \move(1.44338 -5.5)
  \HHSchritt
  \move(3.17543 -4.5)
  \HHSchritt
  \move(0.57735 -7)
  \HHSchritt
  \move(1.44338 -7.5)
  \HHSchritt
  \move(2.3094 -7)
  \HHSchritt
  \move(-1.1547 -6)
  \LUSchritt
  \move(-0.57735 -6)
  \RDSchritt
  \move(0.57735 -6)
  \LDSchritt
  \move(1.1547 -6)
  \RDSchritt
  \move(2.02073 -6.5)
  \RUSchritt
  \move(2.88675 -6)
  \RDSchritt
  \move(3.75278 -6.5)
  \RUSchritt
\esegment
}
\vskip10pt
\centerline{\small c. The corresponding perfect matching of the
  ``dual" graph}
\caption{}
\label{fig:PP10}
\end{figure}

\medskip
The next important method to enumerate plane partitions is actually
a more general method to enumerate {\it perfect matchings of bipartite
  graphs}. A perfect matching of a graph is a collection of pairwise
non-incident edges which cover all vertices of the graph.
To see how plane partitions can indeed be seen as 
perfect matchings of certain graphs, in the triangular grid
inside the bounding hexagon (see Figure~\ref{fig:PP3}.a) 
place a vertex into the centre of each triangle and connect
vertices in adjacent triangles by an edge.
In this way, we obtain a hexagonal graph; see
Figure~\ref{fig:PP10}.a. A unit
rhombus is the union of two adjacent triangles.
Consequently, a rhombus tiling of a given hexagon corresponds
to a collection of edges in the hexagonal graph which are
pairwise non-incident and cover all vertices of the graph;
in other words: a perfect matching of the hexagonal graph.
See Figure~\ref{fig:PP10}.b,c for an example.

Let $G$ be a given graph. If we consider the defining expansion
of the Pfaffian (cf.\ \cite[p.~102]{StemAE}) of the adjacency matrix of
$G$, then each term in this expansion corresponds to a perfect
matching of $G$. Kasteleyn \cite{KastAA,KastAB}
(see \cite{TeslAF} for an excellent exposition) showed that
for planar graphs
there is a consistent way to introduce signs to the $1$'s in
the adjacency matrix so that the Pfaffian of this {\it modified\/}
adjacency matrix (often called {\it Kasteleyn--Percus matrix})
has the property that all terms in its expansion
have the same sign, and hence gives the number of perfect
matchings of~$G$. 
Since in the case of a bipartite graph, this Pfaffian reduces
to a determinant (namely the determinant of the suitably
modified bipartite adjacency matrix of the graph), this
approach has been called {\it``permanent-determinant method"}
by Kuperberg \cite{KupeAA,KupeAG}. He has been the first and
only one to effectively apply this approach for exact
enumeration of plane partitions, providing in particular the first proof
for the enumeration of cyclically symmetric self-complementary
plane partitions (Class~9 in Section~\ref{sec:10}). As a matter of
fact, since the permanent-determinant method produces determinants
(and permanents and Pfaffians) of very large matrices, this approach
is best suited for conceptual insight in the enumeration
of plane partitions (and, more generally, perfect matchings of
bipartite graphs; see e.g.\ the last paragraph in
Section~\ref{sec:endofstory}), whereas less so for proving explicit
exact enumeration formulae. It is interesting to note, though,
that Kuperberg has shown that, for planar
graphs, the non-intersecting lattice path method and the 
permanent-determinant method are equivalent, in a sense explained
in \cite[Sec.~3.3]{KupeAI}. This does not include the situations where the
starting or/and end points of the non-intersecting lattice paths
to be counted vary in given sets (cf.\ the paragraph below
\eqref{eq:MM}); Kuperberg has however also
presented workarounds in \cite{KupeAA} 
to simulate certain such situations within
the permanent-determinant method, and he speculates in 
\cite[end of Sec.~3.3]{KupeAI} that this is always possible (in the
planar case).

\medskip
Extremely useful tools for the enumeration of plane partitions
(and, more generally, perfect matchings)
are moreover Kuo's condensation method \cite{KuoEAA}, which allows
for inductive proofs of (conjectured) enumeration formulae,
Ciucu's matchings factorisation theorem \cite{CiucAB} for
enumerating perfect matchings of graphs with a reflective symmetry, and
Jockusch's theorem \cite{JockAA} for graphs with a rotational symmetry.
As Fulmek \cite{FulmAI} (for the former) and Kuperberg \cite{KupeAG}
(for the latter two) have shown, all three are consequences of the
permanent-determinant method.

\medskip
Finally, among the purely combinatorially methods that have been
applied to the enumeration of plane partition,
Robinson--Schensted--Knuth-like correspondences must be
mentioned; see \cite{GansAC,GansAE,StanAC} and the next section.

\section{Trace generating functions}
\label{sec:trace}

In this section, I present one of the main results from Stanley's 
first article \cite{StanAC} 
on plane partitions.\footnote{It is indeed the chronologically first
article, see Stanley's
publication list {\tt http://www-math.mit.edu/\~{}rstan/pubs/} on his
website. The order of his papers in the References section at the
end of this article follows
this order.}
There, Stanley introduces the notion of ``trace" of
a plane partition. As already said above, the {\it trace} of a plane
partition $\pi$ as in \eqref{eq:PP} is the sum of its parts
along the main diagonal, that is, $\sum_{i\ge1}\pi_{i,i}$.
He is led to this notion by observing that a combination of
a construction of Frobenius \cite[p.~523]{FrobAA} for partitions with
a variation of the
Robinson--Schensted--Knuth correspondence \cite{KnutAA} leads to
an elegant and effortless proof of \eqref{eq:PPall}, and that in
this proof a certain statistic, namely the trace, is preserved, leading
to a refinement of \eqref{eq:PPall} (see \eqref{eq:PPtrace} below).

More precisely, let us consider the plane partition
\begin{equation} \label{eq:PP0} 
\begin{matrix} 
4&3&2&2\\
4&3&1&1\\
2&2&1&1\\
1&1\\
1&1
\end{matrix}
\end{equation}
Each row by itself is a(n ordinary) partition. For example,
the first row represents the partition $4+3+2+2$. A partition
$\la_1+\la_2+\dots+\la_k$
may be represented as a {\it Ferrers diagram}, that is,
as a left-justified diagram of dots, with $\la_i$ dots
in the $i$-th row. Thus, the Ferrers diagram of the
above partition is
\begin{equation} \label{eq:Ferrers} 
\begin{matrix} 
\bullet&\bullet&\bullet&\bullet\\
\bullet&\bullet&\bullet\\
\bullet&\bullet\\
\bullet&\bullet
\end{matrix}
\end{equation}
The partition {\it conjugate} to $\la$ is the partition
$\la'$ corresponding to the transpose of the Ferrers diagram
of $\la$. Thus, the partition conjugate to $4+3+2+2$ is $4+4+2+1$.

By conjugating each row of \eqref{eq:PP0}, we obtain
\begin{equation} \label{eq:PP1}
\begin{matrix} 
4&4&2&1\\
4&2&2&1\\
4&2\\
2\\
2
\end{matrix}
\end{equation}
It should be noted that the trace of the original
plane partition --- in our example in \eqref{eq:PP0} this is $4+3+1=8$
--- equals the number of parts $\pi_{i,j}$ with
$\pi_{i,j}\ge i$ in the plane partition one obtains under this
mapping. Moreover, the total sum of the parts of the plane partition
is preserved under this operation.

Obviously, each {\it column} of the plane partition obtained is
a(n ordinary) partition. We now write each column in (essentially) its
{\it Frobenius notation}. More precisely, out of a partition
$\la_1+\la_2+\dots+\la_k$, we form the pair 
$$
(\la_1,\la_2-1,\la_3-2,\ldots\mid
\la'_1,\la'_2-1,\la'_3-2,\dots),
$$
with $\la'$ denoting the conjugate of $\la$, and only recording
positive numbers.
(In the Ferrers diagram picture, this corresponds to splitting
the partition $\la$ along its main diagonal, ``counting" the main
diagonal ``twice". In Frobenius' original notation, the main
diagonal is {\it not\/} accounted for.) For example, from the
partition $4+3+2+2$
in our running example we obtain the pair $(4,2\mid 4,3)$.
We apply this modified Frobenius notation to each column of
the plane partition we have obtained so far, and we collect
the first components and the second components separately.
Thus, for our running example \eqref{eq:PP1}, we obtain
\begin{equation} \label{eq:PP3}
\begin{matrix} 
4&4&2&1\\
3&1&1\\
2
\end{matrix} 
\quad \quad 
\begin{matrix} 
5&3&2&2\\
4&2&1\\
1
\end{matrix}
\end{equation}
(Here, the pair $(4,3,2\mid 5,4,1)$ corresponds to the first column
in \eqref{eq:PP1}, etc.)
It is easy to see that in this way we obtain a pair $(C_1,C_2)$
of so-called {\it column-strict} plane partitions, which are
plane partitions with the additional condition that parts along
columns are {\it strictly} decreasing. Moreover, $C_1$ and $C_2$
have the same shape. It should be observed that the trace of the
original plane partition equals the number of parts of $C_1$ (or of
$C_2$). It is not true that the total sum of the parts of the original
plane partition equals the sum of all the parts of $C_1$ and $C_2$.
However, what is true is that the total sum of the original plane
partition equals the sum of all parts of $C_1$ and $C_2$ {\it minus}
the {\it number} of parts of $C_1$.

By a variant of the Robinson--Schensted--Knuth correspondence
(see \cite[p.~368]{StanBI}), the pairs $(C_1,C_2)$ correspond to
matrices $(m_{i,j})_{i,j\ge1}$ with a finite number of 
positive integer entries, all other entries being zero.
In this correspondence, the total sum of the parts of $C_1$
and $C_2$ equals $\sum_{i,j\ge1}(i+j)m_{i,j}$ (actually, something
much finer is true), and the number of parts of $C_1$ equals
$\sum_{i,j\ge1}m_{i,j}$. If one puts everything together,
this shows (see \cite[Sec.~3 in combination with Theorem~2.2]{StanAC})
\begin{equation} \label{eq:PPtrace}
\sum_\pi t^{\text{trace}(\pi)}q^{\vert\pi\vert}
=
\prod _{i,j\ge1} ^{}\frac {1} {1-tq^{i+j-1}}. 
\end{equation}
In \cite{StanAC}, we also find the notion of {\it``conjugate trace"},
together with further generating function results for plane partitions
featuring that variation of ``trace."

\medskip
Almost ten years later, Stanley's first student, Emden Gansner, showed
that there is much more to Stanley's concept of ``trace."
He generalised it and defined the 
{\it $i$-trace} of a plane partition $\pi$, denoted by $t\sb
i(\pi)$, as the sum of the parts of the $i$-diagonal of $\pi$,
that is, as $\sum_{\ell\ge1}\pi_{\ell,\ell+i}$. Thus, in this more
general context, Stanley's original trace is the $0$-trace.

Gansner's result on trace generating functions in 
\cite{GansAE}
actually concerns {\it reverse plane partitions}, which are
arrays of non-negative integers as in \eqref{eq:PP} with the
property that entries along rows and along columns are weakly
{\it increasing}. For the formulation of the result, we need
to define the monomial $\mathbf x(\rh;\la)$ for a ``cell"~$\rh$ 
of a partition~$\la$ (a cell corresponds to a bullet in the
Ferrers diagram of~$\la$; see \eqref{eq:Ferrers}). If $\rh$ is located
in the $i$-th row and the $j$-th column of $\la$, then
$$\mathbf x(\rho):=
x_{j-\la'_j}\cdots x_{\la_i-i-1}x_{\la_i-i},
$$
where $\la'$ is again the partition conjugate to $\la$
(see the paragraph after \eqref{eq:Ferrers}). Using this notation,
Gansner's elegant result \cite[Theorem~5.1]{GansAE} is
\begin{equation} \label{eq:Gansner}
\sum_{\pi} 
\prod _{i=1-\ell(\la)} ^{\la_1-1}x_i^{t_i(\pi)}
=
\prod _{\rho\in\la} ^{}\frac {1} {1-\mathbf x(\rho)},
\end{equation}
where the sum is over all reverse plane partitions $\pi$ of shape~$\la$.
The proof is bijective; it exploits hidden properties of a
correspondence of Hillman and Grassl \cite{HiGrAA} that Gansner
lays open. He continued that work in \cite{GansAE}, where he
considers trace generating functions for (ordinary) plane partitions.
As it turns out, the trace generating function (with all traces present)
for plane partitions of a given shape does not evaluate to an
equally elegant product formula as for reverse plane partitions
(except for rectangular shapes, which however follows trivially from the
result \eqref{eq:Gansner} for reverse plane partitions by a
$180^\circ$ degree rotation), one has to be content with formulae
given in terms of certain sums. In order to derive these results,
Gansner uses Robinson--Schensted--Knuth-like correspondences due to Burge
\cite{BurgAA}.
It is shown in \cite[Theorems~7.6, 7.7]{KratAM} that, for generating
functions for plane partitions which only take into account the size
and $0$-trace of plane partitions, there are closed form products
available for certain shapes.

Finally, we mention that Proctor \cite{ProcAB} has derived
several results for {\it alternating trace} generating functions
by relating these to characters of symplectic groups.

\section{The 10 symmetry classes of plane partitions}
\label{sec:10}

The purpose of this section is to summarise the programme
of enumeration of symmetry classes of plane partitions,
the history of which was reviewed in Section~\ref{sec:Stan}.
I start by defining the symmetry operations which give rise to
10 symmetry classes of plane partitions. Then, for each
symmetry class, I state the corresponding enumeration formula,
and report who had proved that formula, sometimes with pointers
to further work. 

The first operation is {\it reflection}. To reflect
a plane partition $\pi$ as in \eqref{eq:PP}
means to reflect it along the main diagonal, that is,
to map $\pi=(\pi_{i,j})_{i,j}$ to $\pi=(\pi_{j,i})_{i,j}$.
In the representation of a plane partition
as a rhombus tiling as in Figure~\ref{fig:PP3}.b, reflection 
means reflection along the vertical symmetry
axis of the hexagon (if there is one). 

{\it Rotation} is the second operation.
Given a plane partition $\pi$, viewed as a pile of cubes as
in Figure~\ref{fig:PP2}, to rotate it means to rotate it
by $120^\circ$ with rotation axis 
$\{(t,t,t):-\infty<t<\infty\}$ (this is the rotation
which leaves the origin of the coordinate
system invariant and maps the coordinate axes to each other;
for our purposes it is irrelevant whether we consider a
``left" or ``right" rotation).
In the representation of a plane partition
as a rhombus tiling as in Figure~\ref{fig:PP3}.b, this rotation
corresponds to a rotation of the
hexagon (together with the tiling) by $120^\circ$.

Intuitively, the {\it complement\/} of a plane partition in a given
$a\times b\times c$ box is what the name says, namely the
pile of cubes which remains (suitably rotated and reflected) 
after the plane partition is removed from the box. To be
precise, if the cubes of a plane partition $\pi$ contained in
an $a\times b\times c$ box are coordinatised in the
obvious way by $(i,j,k)$ with $1\le i\le a$, $1\le j\le b$, and $1\le
k\le c$, then the complement $\pi^c$ of $\pi$ is defined by
$$
\pi^c=\{(a+1-i,b+1-j,c+1-k):(i,j,k)\notin\pi\}.
$$
In the representation of a plane partition
as a rhombus tiling as in Figure~\ref{fig:PP3}.b, 
to take the complement of a plane partition means to rotate the
hexagon (together with the tiling) by $180^\circ$.

If one combines the three operations --- reflection, rotation,
complementation --- in all possible ways, then there result
ten {\it symmetry classes} of plane partitions. Remarkably, in
all ten cases, there exist closed form product formulae for the
number of plane partitions contained in a given box, respectively
for certain generating functions. The first complete
presentation of the corresponding results and conjectures 
was given by Stanley in \cite{StanAI}, together
with the state of the art at the time, and the proof of
the conjectured formulae for Class~5 (see also Section~\ref{sec:SCPP}).

While reading the descriptions of the various symmetry
classes, Figure~\ref{fig:PP4} may be helpful. 
It shows a plane partition which has {\it all\/} possible
symmetries. Thus, it belongs to all ten classes.

\begin{figure}
\begin{center}
\setlength{\unitlength}{.10cm}

\end{center}
\caption{A totally symmetric self-complementary plane partition\newline
--- as a pile of cubes, and as a rhombus tiling of a hexagon}
\label{fig:PP4}
\end{figure}

In order to have a convenient
notation to write down some of the formulae, we use the symbol
$N_d(a,b,c)$ for the number of all plane partitions in symmetry
class~$d$ which are contained in an $a\times b\times c$ box. 

\medskip
{\sc Class 1:} {\it Unrestricted Plane Partitions.}
As we already said, MacMahon proved
(see Theorem~\ref{thm:MM}) that the generating function for all plane
partitions contained in an $a\times b\times c$ box is given by
\begin{equation} \label{eq:C1} 
\sum_{\pi}q^{\vert\pi\vert}=
\prod _{i=1} ^{a}
\prod _{j=1} ^{b}
\prod _{k=1} ^{c}
\frac {1-q^{i+j+k-1}} {1-q^{i+j+k-2}},
\end{equation}
where the sum is over all plane partitions $\pi$ contained in an
$a\times b\times c$ box.

\medskip
{\sc Class 2:} {\it Symmetric Plane Partitions.}
A plane partition $\pi$ as in \eqref{eq:PP}
is called {\it symmetric} if it is invariant under
reflection along the main diagonal.
In the representation of a plane partition
as a rhombus tiling as in Figure~\ref{fig:PP3}.b, to be symmetric
means to be invariant under reflection along the vertical symmetry
axis of the hexagon 
(two of the side lengths of the hexagon must be equal). 

For this class, two different weights lead to generating functions
which have closed form product formulae. The first weight is the usual
``size" of a plane partition, and it leads to MacMahon's conjecture
\cite{MacMAC} that 
\begin{equation} \label{eq:C2a}
\sum_\pi q^{\vert\pi\vert}=
\prod _{i=1} ^{a}\frac {1-q^{c+2i-1}} {1-q^{2i-1}}
\prod _{1\le i<j\le a} ^{}\frac {1-q^{2(c+i+j-1)}} {1-q^{2(i+j-1)}} ,
\end{equation}
where the sum is over all symmetric plane partitions $\pi$ that
are contained in an $a\times a\times c$ box.
MacMahon's conjecture was first proved by
Andrews \cite{AndrAK} and Macdonald \cite[Ex.~16 and 17,
pp.~83--85]{MacdAC}.
Since then, several other proofs and refinements were
given, see \cite[Prop.~7.3]{ProcAD}, \cite[Theorem~1, Case~BYI]{ProcAB},
\cite[Sec.~5]{KupeAH}, \cite[Cor.~12]{KratAP}, \cite[Theorem~5]{KratBC}.

The second weight is roughly ``half" the size of a plane partition,
and it leads to a conjecture of Bender and Knuth \cite[p.~50]{BeKnAA}.
Given a plane partition $\pi$ as in \eqref{eq:PP}, 
the weight $\vert\pi\vert_0$ of $\pi$
is defined as the sum $\sum_{1\le j\le i} \pi_{i,j}$. Then
\begin{equation} \label{eq:C2b}
\sum_\pi q^{\vert\pi\vert_0}=
\prod _{1\le i\le j\le a} ^{}\frac {1-q^{c+i+j-1}} {1-q^{i+j-1}},
\end{equation}
where the sum is over all symmetric plane partitions $\pi$ that
are contained in an $a\times a\times c$ box.
Bender and Knuth's conjecture was first proved by
Gordon \cite{GordAC} (as reported in \cite[Prop.~16.1]{StanAA}), 
but published
only much later. The first published proofs are due to Andrews
\cite{AndrAJ} and Macdonald \cite[Ex.~19, p.~86]{MacdAC}.
Since then, several other proofs and refinements were
given, see \cite[Prop.~7.2]{ProcAD}, \cite[Theorem~1,
  Case~BYH]{ProcAB}, 
\cite[Theorem~21]{KratAP}, \cite[Theorem~6]{KratBC}, \cite{FiscAH}.

\medskip
{\sc Class 3:} {\it Cyclically Symmetric Plane Partitions.}
A plane partition $\pi$
is called {\it cyclically symmetric} if, viewed as a pile of cubes as
in Figure~\ref{fig:PP2}, it is invariant under 
rotation by $120^\circ$ as described above. 
In the representation of a plane partition
as a rhombus tiling as in Figure~\ref{fig:PP3}.b, to be cyclically
symmetric means to be invariant under rotation of the
hexagon (together with the tiling) 
by $120^\circ$ (all sides of the hexagon must have the
same length).

For this class, Macdonald \cite[Ex.~18, p.85]{MacdAC} conjectured
\begin{equation} \label{eq:C3}
\sum_\pi q^{\vert\pi\vert}=
\prod _{i=1} ^{a}\frac {1-q^{3i-1}} {1-q^{3i-2}}
\prod _{1\le i<j\le a} ^{}\frac {1-q^{3(2i+j-1)}}
      {1-q^{3(2i+j-2)}}
\prod _{1\le i<j,k\le a} ^{}\frac {1-q^{3(i+j+k-1)}}
      {1-q^{3(i+j+k-2)}}, 
\end{equation}
where the sum is over all cyclically symmetric plane partitions $\pi$ that
are contained in an $a\times a\times a$ box.

Andrews \cite[Theorem~4]{AndrAO} had found a determinant for the generating
function in question; in retrospective this determinant can be
explained by non-intersecting lattice paths (which have been 
discussed in Section~\ref{sec:methods}).
The conjecture was proved by Mills, Robbins and Rumsey \cite{MiRRAA}
by evaluating this determinant
using some beautiful linear algebra, and it is still the only proof.

\medskip
{\sc Class 4:} {\it Totally Symmetric Plane Partitions.}
A plane partition is called {\it totally symmetric} if it is at the
same time symmetric and cyclically symmetric.

For this class, Andrews and Robbins (see \cite{AndrZZ}) conjectured
\begin{equation} \label{eq:C4}
\sum_\pi q^{\vert\pi\vert_0}=
\prod _{1\le i\le j\le k\le a} ^{}
\frac {1-q^{i+j+k-1}} {1-q^{i+j+k-2}},
\end{equation}
where the sum is over all totally symmetric plane partitions $\pi$ that
are contained in an $a\times a\times a$ box.

Okada \cite[Theorem~5]{OkadAA} found a determinant for the generating
function in question. However, for a long time nobody knew how to
evaluate this determinant. In the $q=1$
special case, Stembridge \cite{StemAG} was able
to relate a slightly different determinant to the enumeration
of cyclically symmetric plane partitions, thereby solving the problem of 
plain enumeration of totally symmetric plane partitions.
Unfortunately, it seems that this very conceptual approach does not extend to the
$q$-case. This was the longest standing conjecture of all plane
partition conjectures.
Finally, Koutschan, Kauers and Zeilberger \cite{KoKZAC}
succeeded to evaluate Okada's determinant, by a (heavily)
computer-assisted argument. 

\medskip
{\sc Class 5:} {\it Self-Complementary Plane Partitions.}
A plane partition is called {\it self-comple\-mentary} if it is equal to
its complement.
In the representation of a plane partition
as a rhombus tiling as in Figure~\ref{fig:PP3}.b, to be self-complementary
means to be invariant under rotation of the
hexagon (together with the tiling) 
by $180^\circ$. In other words, to be self-complementary in the
rhombus tiling picture means to be centrally symmetric.

For this class, Robbins and Stanley independently conjectured that
{\refstepcounter{equation}\label{eq:C5}}
\alphaeqn
\begin{align} 
\label{eq:C5a}
N_5(2a,2b,2c)&=N_1(a,b,c)^2,\\
\label{eq:C5b}
N_5(2a+1,2b,2c)&=N_1(a,b,c)N_1(a+1,b,c),\\
\label{eq:C5c}
N_5(2a+1,2b+1,2c)&=N_1(a+1,b,c)N_1(a,b+1,c).
\end{align}
\reseteqn
As mentioned above, Stanley proved this conjecture in
\cite[Sec.~3]{StanAI}, see Section~\ref{sec:SCPP}.
For a different proof see \cite[Sec.~4]{KupeAH} and \cite[Ex.~4.2]{StemAI}.

\medskip
{\sc Class 6:} {\it Transpose-Complementary Plane Partitions.}
A plane partition is called {\it transpose-comple\-mentary} if it is equal to
the reflection of its complement.
In the representation of a plane partition
as a rhombus tiling as in Figure~\ref{fig:PP3}.b, to be transpose-complementary
means to be invariant under reflection of the
hexagon (together with the tiling) in its horizontal symmetry axis
(two of the side lengths of the hexagon must be equal).

The number of all transpose-complementary plane partitions
contained in an $a\times a\times c$ box is equal to
\begin{equation} \label{eq:C6}
\binom {c+a-1}{a-1}
\prod _{1\le i\le j\le a-2} ^{}\frac {2c+i+j+1} {i+j+1}.
\end{equation}
This was proved by Proctor \cite[paragraph above Cor.~1,
Cases `CG']{ProcAB} by (in essence) relating
these plane partitions to the symplectic tableaux
of King and El-Sharkaway \cite{KiElAA}.
For a different proof see \cite[Sec.~6]{KupeAH}.

\medskip
{\sc Class 7:} {\it Symmetric Self-Complementary Plane Partitions.}
A plane partition is called {\it symmetric self-comple\-mentary} 
if it is equal to its complement and invariant under reflection in
the main diagonal.
In the representation of a plane partition
as a rhombus tiling as in Figure~\ref{fig:PP3}.b, to be 
symmetric self-complementary
means to be invariant under reflection of the
hexagon (together with the tiling) in its vertical symmetry axis
and under reflection in its horizontal symmetry axis
(all side lengths of the hexagon must be equal).
For this class, we have
{\refstepcounter{equation}\label{eq:C7}}
\alphaeqn
\begin{align} 
\label{eq:C7a}
N_7(2a,2a,2c)&=N_1(a,a,c),\\
\label{eq:C7b}
N_7(2a+1,2a+1,2c)&=N_1(a+1,a,c).
\end{align}
\reseteqn
This was proved by Proctor \cite{ProcAE} by again using the
representation theory of symplectic groups.
For a different proof see \cite[Ex.~4.3]{StemAI}.

\medskip
{\sc Class 8:} {\it Cyclically Symmetric Transpose-Complementary 
Plane Partitions.}
A plane partition is called {\it cyclically symmetric transpose-comple\-mentary} 
if it is equal to the transpose of its complement and invariant under
rotation  
by $120^\circ$.
In the representation of a plane partition
as a rhombus tiling as in Figure~\ref{fig:PP3}.b, to be 
cyclically symmetric transpose-complementary
means to be invariant under rotation of the
hexagon (together with the tiling) by $120^\circ$
and under reflection in its horizontal symmetry axis
(all side lengths of the hexagon must be equal).

The number of all cyclically symmetric transpose-complementary plane partitions
contained in a $2a\times 2a\times 2a$ box is equal to
\begin{equation} \label{eq:C8}
\prod _{i=0} ^{a-1}\frac {(3i+1)!\,(6i)!\,(2i)!} 
{(4i+1)!\,(4i)!}. 
\end{equation}
This was proved by Mills, Robbins and Rumsey \cite{MiRRAD}
by an interesting determinant evaluation.
For a different proof, using Ciucu's matchings factorisation theorem
\cite{CiucAB}, see \cite[Sec.~4]{CiKrAC}.

\medskip
{\sc Class 9:} {\it Cyclically Symmetric Self-Complementary 
Plane Partitions.}
A plane partition is called {\it cyclically symmetric self-comple\-mentary} 
if it is equal to its complement and invariant under rotation
by $120^\circ$.
In the representation of a plane partition
as a rhombus tiling as in Figure~\ref{fig:PP3}.b, to be 
cyclically symmetric self-complementary
means to be invariant under rotation of the
hexagon (together with the tiling) by $60^\circ$
(all side lengths of the hexagon must be equal).

The number of all cyclically symmetric self-complementary plane partitions
contained in a $2a\times 2a\times 2a$ box is equal to
\begin{equation} \label{eq:C9}
\prod _{i=0} ^{a-1}\frac {(3i+1)!^2} {(a+i)!^2}.
\end{equation}
This was first proved by Kuperberg \cite{KupeAA} by using Kasteleyn's
method (cf.\ \cite{KastAA,KastAB}).
Ciucu \cite{CiucAI} found a combinatorial explanation of the
factorisation 
\begin{equation} \label{eq:C9C10}
N_9(2a,2a,2a)=N_{10}(2a,2a,2a)^2
\end{equation}
that results from a comparison of \eqref{eq:C9} and \eqref{eq:C10},
based on his matching factorisation theorem
\cite{CiucAB}.
A direct proof of \eqref{eq:C9}, using Ciucu's matchings factorisation theorem
\cite{CiucAB}, can be found in \cite[Sec.~5]{CiKrAC}.

\medskip
{\sc Class 10:} {\it Totally Symmetric Self-Complementary 
Plane Partitions.}
A plane partition is called {\it totally symmetric self-comple\-mentary} 
if it is equal to its complement, is invariant under rotation by
$120^\circ$, and also under reflection in
the main diagonal.
In the representation of a plane partition
as a rhombus tiling as in Figure~\ref{fig:PP3}.b, to be 
totally symmetric self-complementary means to be invariant under all
symmetry operations of the hexagon, that is, under reflection of the
hexagon (together with the tiling) in its horizontal symmetry axis
and under rotation by $60^\circ$
(all side lengths of the hexagon must be equal).

The number of all totally symmetric self-complementary plane partitions
contained in a $2a\times 2a\times 2a$ box is equal to
\begin{equation} \label{eq:C10}
\prod _{i=0} ^{a-1}\frac {(3i+1)!} {(a+i)!}.
\end{equation}
Stembridge \cite[Theorem~8.3]{StemAE} found a determinant for this number
using non-intersecting lattice paths,
and Andrews \cite{AndrAW} succeeded to evaluate this determinant.
See \cite{AnBuAA} for a much simpler proof of that determinant
evaluation, and
\cite{KratBD} for a generalisation containing an additional parameter.

\section{Stanley's proof of the formula for self-complementary
plane partitions}
\label{sec:SCPP}

Here is a ``modern" version of Stanley's proof of \eqref{eq:C5},
``modern" in the sense that we argue using the rhombus tiling
point of view of plane partitions. We want to count all
self-complementary plane partition contained in an $a\times b\times c$
box, or, equivalently, all rhombus tilings of a hexagon with side
lengths $a,b,c,a,b,c$ which are invariant under rotation by
$180^\circ$. There are in fact three cases
to consider: all three of $a,b,c$ are even, two of $a,b,c$ are even,
or only one of $a,b,c$ is even. Clearly, if all of $a,b,c$ are odd,
then there does not exist such a self-complementary plane partition
(the plane partition together with its complement have to fill the
complete $a\times b\times c$ box, which consequently must have
even volume). 

\begin{figure}[h]
\centertexdraw{
  \drawdim truecm  \linewd.02
\move(0 0)
\bsegment
  \drawdim truecm  \linewd.02
  \linewd.12
  \move(0 0)
  \RhombusA \RhombusA \RhombusB 
  \RhombusA \RhombusB \RhombusA \RhombusB \RhombusA \RhombusB
  \RhombusA
  \move (-0.866025403784439 -.5)
  \RhombusA \RhombusB \RhombusA \RhombusB 
  \RhombusA \RhombusB \RhombusA \RhombusB 
  \RhombusA \RhombusA
  \move (-0.866025403784439 -.5)
  \RhombusC 
  \move (-0.866025403784439 -1.5)
  \RhombusC \RhombusC
  \move (-0.866025403784439 -2.5)
  \RhombusC \RhombusC \RhombusC 
  \move (-0.866025403784439 -3.5)
  \RhombusC \RhombusC \RhombusC \RhombusC 
  \move (2.598076211353316 -.5)
  \RhombusC \RhombusC \RhombusC \RhombusC 
  \move (3.4641 -2)
  \RhombusC \RhombusC \RhombusC 
  \move (4.330125 -3.5)
  \RhombusC \RhombusC 
  \move (5.19615 -5)
  \RhombusC 
\esegment
\move(8 0)
\bsegment
  \drawdim truecm  \linewd.02
  \linewd.12
  \move(3.4641 -4)
  \RhombusB \RhombusA \RhombusB
  \RhombusA
  \move (0.866025403784439 -2.5)
  \RhombusB 
  \RhombusA \RhombusB \RhombusA \RhombusB 
  \RhombusA \RhombusA
  \move (-0.866025403784439 -1.5)
  \RhombusC \RhombusC
  \move (-0.866025403784439 -2.5)
  \RhombusC \RhombusC \RhombusC 
  \move (-0.866025403784439 -3.5)
  \RhombusC \RhombusC \RhombusC \RhombusC 
  \move (2.598076211353316 -3.5)
  \RhombusC 
  \move (5.19615 -5)
  \RhombusC 
\esegment}
\vskip15pt
\centerline{a. A self-complementary plane partition\hskip2cm
b. Its lower ``half"\hskip1cm}
\caption{}
\label{fig:PP5}
\end{figure}

We focus on the case where all of $a,b,c$ are even, say $a=2a_1$,
$b=2b_1$, $c=2c_1$, all other cases being
similar. Figure~\ref{fig:PP5}.a
shows an example of a rhombus tiling contained in a $2\times 6\times
4$ box which is invariant under rotation by $180^\circ$. Obviously, 
one half of the tiling determines the rest. Therefore we may
concentrate on just the lower half, see Figure~\ref{fig:PP5}.b. 
It should be noted
that there are some rhombi of the original tiling which stick out of
the half-hexagon along the cut. In fact, there are exactly
$\min\{a,c\}$ such rhombi. Moreover, since the original
tiling was centrally symmetric, these rhombi must be arranged
symmetrically along the cut. Now, as in Section~\ref{sec:methods}
(see Figure~\ref{fig:PP8}), we mark the
midpoints of the edges along the south-west side of the half-hexagon,
and we start paths there, where the individual steps of the paths
always connect midpoints of opposite sides of rhombi. See
Figure~\ref{fig:PP6}
for the result in our running example. We obtain a collection of paths
which connect the midpoints of the edges along the south-west side
with the midpoints of parallel edges along the cut. Clearly, the
paths are {\it non-intersecting}, meaning that no two paths have any
points in common. Each path consists of $\frac {1} {2}(a+c)$ steps.

\begin{figure}[h]
\centertexdraw{
  \drawdim truecm  \linewd.02
\move(0 0)
\bsegment
  \drawdim truecm  \linewd.02
  \linewd.12
  \move(3.4641 -4)
  \RhombusB \RhombusA \RhombusB
  \RhombusA
  \move (0.866025403784439 -2.5)
  \RhombusB 
  \RhombusA \RhombusB \RhombusA \RhombusB 
  \RhombusA \RhombusA
  \move (-0.866025403784439 -1.5)
  \RhombusC \RhombusC
  \move (-0.866025403784439 -2.5)
  \RhombusC \RhombusC \RhombusC 
  \move (-0.866025403784439 -3.5)
  \RhombusC \RhombusC \RhombusC \RhombusC 
  \move (2.598076211353316 -3.5)
  \RhombusC 
  \move (5.19615 -5)
  \RhombusC 
\linewd.08
\move(-.433012 -4.75)
\vdSchritt \vdSchritt \vdSchritt 
\move(.433012 -5.25)
\vdSchritt \vdSchritt \vdSchritt 
\move(1.29904 -5.75)
\vdSchritt \vdSchritt \hdSchritt 
\move(2.165065 -6.25)
\vdSchritt \hdSchritt \vdSchritt 
\move(3.031090 -6.75)
\hdSchritt \vdSchritt \hdSchritt 
\move(3.897115 -7.25)
\hdSchritt \hdSchritt \vdSchritt 
\ringerl(-.433012 -4.75) 
\ringerl(.433012 -5.25)
\ringerl(1.29904 -5.75)
\ringerl(2.165065 -6.25)
\ringerl(3.031090 -6.75)
\ringerl(3.897115 -7.25)
\ringerl(-.433012 -1.75) 
\ringerl(.433012 -2.25)
\ringerl(2.165065 -3.25)
\ringerl(3.031090 -3.75)
\ringerl(4.763140 -4.75)
\ringerl(5.629165 -5.25)
\esegment
\move(8 0)
\bsegment
  \drawdim truecm  \linewd.02
  \linewd.12
\linewd.08
\move(-.433012 -4.75)
\vdSchritt \vdSchritt \vdSchritt 
\move(.433012 -5.25)
\vdSchritt \vdSchritt \vdSchritt 
\move(1.29904 -5.75)
\vdSchritt \vdSchritt \hdSchritt 
\move(2.165065 -6.25)
\vdSchritt \hdSchritt \vdSchritt 
\move(3.031090 -6.75)
\hdSchritt \vdSchritt \hdSchritt 
\move(3.897115 -7.25)
\hdSchritt \hdSchritt \vdSchritt 
\ringerl(-.433012 -4.75) 
\ringerl(-.433012 -3.75) 
\ringerl(-.433012 -2.75) 
\ringerl(.433012 -5.25)
\ringerl(.433012 -4.25)
\ringerl(.433012 -3.25)
\ringerl(1.29904 -5.75)
\ringerl(1.29904 -4.75)
\ringerl(1.29904 -3.75)
\ringerl(2.165065 -6.25)
\ringerl(2.165065 -5.25)
\ringerl(3.031090 -4.75)
\ringerl(3.031090 -6.75)
\ringerl(3.897115 -5.25)
\ringerl(3.897115 -6.25)
\ringerl(3.897115 -7.25)
\ringerl(4.763140 -6.75)
\ringerl(5.629165 -6.25)
\ringerl(-.433012 -1.75) 
\ringerl(.433012 -2.25)
\ringerl(2.165065 -3.25)
\ringerl(3.031090 -3.75)
\ringerl(4.763140 -4.75)
\ringerl(5.629165 -5.25)
\htext(4.297115 -7.5){1}
\htext(5.163140 -7){2}
\htext(3.431090 -7){1}
\htext(4.297115 -5.5){3}
\htext(2.565065 -5.5){2}
\htext(1.79904 -4){3}
\esegment}
\vskip10pt
\centerline{a. Half of a self-complementary plane partition\hskip1cm
b. Non-intersecting lattice paths}
\caption{}
\label{fig:PP6}
\end{figure}

\begin{figure}[h]
$$
\begin{matrix} 
1&1&2&3\\
2&3
\end{matrix}
$$
\caption{The corresponding semistandard tableau}
\label{fig:PP7}
\end{figure}

The next step is to translate these sets of paths into so-called semistandard
tableaux. In order to do this, we label ``right steps" of paths
by their distance to the south-west side of the half-hexagon. More
precisely, a right step in a path is labelled by $i$ if it is the $i$-th
step of this path. Up-steps remain unlabelled. See
Figure~\ref{fig:PP6}.b for the
labelling in our running example. Now we read the labels of the paths
in order, and put the labels of the bottom-most path into the first
column, the labels of the path next to it into the second column, etc.,
of a tableau. Figure~\ref{fig:PP7} shows the tableau which we obtain for the
tiling in Figure~\ref{fig:PP5}.a. By construction, the entries of the tableau 
are strictly increasing along columns. Furthermore, it is not
difficult to see that entries along rows are weakly increasing.
As we said in Section~\ref{sec:methods},
these two properties define a {\it semistandard tableau}. The shape of the
tableau we obtain (where shape is defined in the same way as for a
plane partition) is determined by the position of the rhombi which
stick out of the cut line of the half-hexagon. Thus, it is not
unique, but on the other hand it cannot be arbitrary since it inherits 
the property of these rhombi to be arranged symmetrically along the
cut. To be precise, the only shapes which are possible are the
ones of the form
\begin{equation} \label{eq:shapes}
(b_1+\de_1,b_1+\de_2,\dots,b_1+\de_{a_1},b_1-\de_{a_1},
\dots,b_1-\de_{2},b_1-\de_1) ,
\end{equation}
where $(\de_1,\de_2,\dots,\de_{a_1})$ ranges over all sequences
with $b_1\ge\de_1\ge\de_2\ge\dots\ge\de_{a_1}\ge0$.

If we summarise the discussion so far, we see that the number of
self-complementary plane partitions contained in a $2a_1\times
2b_1\times 2c_1$ box is equal to
\begin{equation} \label{eq:sum1}
{\sum _{\la} ^{}}{}^{}
\left\vert\mathcal T_\la\left(\tfrac {1} {2}(a+c)\right)\right\vert, 
\end{equation}
where $\mathcal T_\la(m)$ denotes the set of semistandard tableaux of shape
$\la$ with entries between $1$ and $m$, and the sum is over all
partitions $\la$ in \eqref{eq:shapes}. By \eqref{eq:Schur}, we see
that $\vert \mathcal T_\la(m)\vert$ equals a specialised Schur
function. More precisely, we have
\begin{equation} \label{eq:Ts} 
\vert\mathcal T_\la(m)\vert=s_\la(1,1,\dots,1),
\end{equation}
with $m$ occurrences of~$1$. Thus, the sum in \eqref{eq:sum1} is a sum
of specialised Schur functions. The crucial observation of Stanley is
that in fact
\begin{equation} \label{eq:s^2} 
{\sum _{\la} ^{}}{}^{}
s_\la(x_1,x_2,\dots)=
s_{(b_1,b_1,\dots,b_1)}(x_1,x_2,\dots)^2,
\end{equation}
with $a_1$ repetitions of $b_1$,
where the sum is over all partitions $\la$ in \eqref{eq:shapes}.
Once this is observed, it is readily verified by means of the
Littlewood--Richardson rule (cf.\ \cite[Ch.~1, Sec.~9]{MacdAC}).
Thus, in view of \eqref{eq:sum1}--\eqref{eq:s^2}, all we
have to do is to set $x_i=1$ for $i=1,2,\dots,\hbox{$\frac {1} {2}(a+c)$}$ in
\eqref{eq:s^2} and evaluate the specialised Schur function on the
right-hand side. How the latter is done, was reviewed in
the paragraph following \eqref{eq:MMSchur}. (In the result, one
has to perform the limit $q\to1$.) 
Thus we have proved \eqref{eq:C5a}.

\section{All problems solved?}
\label{sec:endofstory}

Is this the ``end of the story"? Or is there still something left?
Yes, as reported in Section~\ref{sec:10}, all conjectures on the enumeration of
symmetry classes of plane partitions have now (finally) been
established. However, although significant advances have been made
since Richard Stanley studied plane partitions,
one cannot claim that we have a good understanding of the formulae.
Namely, an immediate question which poses itself is whether there
are insightful explanations why, for all ten symmetry classes, there are
closed form product formulae for the number or the generating function
of plane partitions in the class which are contained in a given box?
Representation theory, in the form of the observation that certain
plane partitions index bases of representations of
classical groups and the fact that there are closed form product
formulae for the dimensions of such representations 
provide explanations for several symmetry classes,
namely for Classes~1, 2, 3 (only for $q=1$), 5, 6, 7, 8
(see \cite{KupeAH,KupeAF,ProcAB} and the references therein). 
Not only is it not one idea that works for all these classes but
several, each of which working for a subset,
for the remaining cases, no such explanation is available. The ``worst" case
is certainly the class of totally symmetric plane partitions,
where the proof consists of a heavily computer-assisted
verification of a (complicated) determinant evaluation.

A combinatorialist would dream of a combinatorial (in the best case:
bijective) argument which would explain the product formulae. However, 
we are very far off a realisation of that dream. Only in the simplest
case, namely the case of MacMahon's theorem for the generating
function for {\it all\/} plane partitions in a given box 
--- presented here in\break Theorem~\ref{thm:MM} ---,
there exists a bijective proof (see \cite{KratBK}). In all other
cases, to find a bijective argument is a wide open problem. 
As mentioned earlier,
there exists at least a combinatorial explanation of the factorisation
\eqref{eq:C9}
occurring for cyclically symmetric self-complementary plane
partitions, see \cite{CiucAI}.

However, the greatest, still unsolved, mystery 
concerns the question what plane
partitions have to do with {\it alternating sign matrices}
(see Section~\ref{sec:Stan} for their definition).
This question was first raised in \cite{StanAH}, and there is
still no answer to it.
Mills, Robbins and Rumsey \cite[Conj.~1]{MiRRAB} conjectured and Zeilberger
\cite{ZeilBD} (see also \cite{KupeAG}) proved that the number
of $n\times n$ alternating sign matrices is given by \eqref{eq:C10}
(with $a$ replaced by $n$),
which also counts totally symmetric self-complementary plane
partitions contained in a $2n\times 2n\times 2n$ box. Moreover,
as Andrews \cite{AndrAO} had shown, the same numbers also arise
as the numbers of ``size~$n$" {\it descending plane partitions}
(which I will not define here). Can these be coincidences?
Certainly not. However, so far nobody has found any conceptual
explanation. 

As a matter of fact, already Mills, Robbins and Rumsey looked deeper
into this mysterious triple occurrence of the same fascinating
number sequence. They attempted to find parametric refinements
and variations of these ``coincidences." Their idea was that,
if the numerical equality took place even at a finer level,
then this may be the guide to the sought-for conceptual
explanation(s) and, in the best case, to bijections between these
objects. The result of their search for refinements was 
several conjectures, one \cite[Conj.~3]{MiRRAB}
connecting refined enumerations
of descending plane partitions and of alternating sign matrices
with each other, and several \cite[Conjectures~2--7']{MiRRAC}
connecting refined enumerations
of totally symmetric self-complementary 
plane partitions and of alternating sign matrices
with each other. The first conjecture has been recently
demonstrated by Behrend, Di~Francesco and Zinn--Justin
\cite{BeDZAA} (and even further refined in \cite{BeDZAB}) 
in a remarkable combination of the new ideas
developed by the second and third author in their
(so far, unsuccessful) attempt to prove the so-called 
Razumov--Stroganov conjecture (cf.\ \cite{CaSpAA}) with
a skillful determinant calculation. As this description
indicates, this proof is again far from being bijective
or illuminating as to why the numerical equality occurs.

The second group of conjectures is still wide open, although
Zeilberger's magnum opus \cite{ZeilBD} actually proves
a weak version of \cite[Conj.~7']{MiRRAC}. 
Again, the proof is highly non-bijective
and non-illuminating. In 1996, the author of these notes found a
generalisation of this second conjecture of Mills, Robbins
and Rumsey. It is also still open, and in particular 
it has so far not led to the sought-for bijection between
alternating sign matrices and totally symmetric self-complementary
plane partitions. Nevertheless, since this conjecture has not
yet been published anywhere, I use the opportunity to present
it here.

\begin{figure}[h]
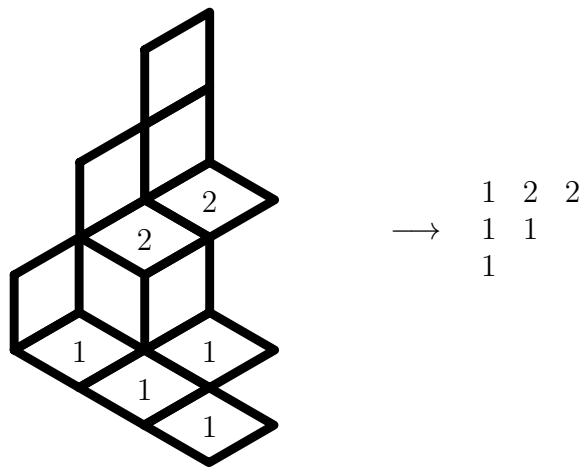

\centertexdraw{
  \drawdim truecm  \linewd.02
\move(8 0)
\bsegment
  \drawdim truecm  \linewd.02
  \linewd.12
  \move(0 0)
  \RhombusB \RhombusB \RhombusA 
  \move (-0.866025403784439 -1.5)
  \RhombusB \RhombusA \RhombusB \RhombusA 
  \move (-1.732050807568877 -3)
  \RhombusB \RhombusA \RhombusA \RhombusA
\htext(-0.12 -2.65){2}
\htext( 0.746025 -2.15){2}
\htext(-0.12 -4.65){1}
\htext( 0.746025 -4.15){1}
\htext( 0.746025 -5.15){1}
\htext(-0.986025 -4.15){1}
\htext( 0.746025 -3.15){$\hskip2.5cm\longrightarrow\quad 
\begin{matrix} 1&2&2\\ 1&1\\1
\end{matrix}$\hskip5cm}
\move (7 0)
\esegment
}
\caption{A totally symmetric self-complementary plane partition 
turned into a triangular array}
\label{fig:Magog}
\end{figure}

The starting point for \cite[Conjecture~7']{MiRRAC} and its
announced generalisation is the fact that both totally symmetric
self-complementary plane partitions and alternating sign matrices
can be encoded in terms of certain triangular arrays.
Namely, a totally symmetric self-complementary plane partition
(see Figure~\ref{fig:PP4}) contained in a $2n\times 2n\times 2n$ box 
is determined by what we see in
just one twelveth of the hexagon, everything else is forced
by symmetry. See Figure~\ref{fig:Magog} for one twelveth of
the plane partition in Figure~\ref{fig:PP4}. The top-most
point in Figure~\ref{fig:Magog} corresponds to the centre in
Figure~\ref{fig:PP4}. Although this is
just one twelveth, we still view this in three dimensions as
earlier. For the ``flat" rhombi, we record their heights in 
relation to the ``floor" (in the three-dimensional
picture) which we assume at height~1;
see the numbers in the left half of 
Figure~\ref{fig:Magog}. These numbers are then arranged in
a triangular array as in the right half of Figure~\ref{fig:Magog}. 
The entries in the triangular array are positive integers. they
have the property that they
are weakly increasing along rows and weakly decreasing 
along columns, and the entries in column~$i$ are bounded
above by~$i$, $i=1,2,\dots,n$.
It is easy to see that the above defines a bijection.

The ``Magog trapezoids" (this is terminology
stolen from \cite{ZeilBD}) below generalise these triangular
arrays. In view of the above correspondence, totally symmetric
self-complementary plane partitions correspond to
$(0,n,n)$-Magog trapezoids.

\begin{definition} \label{def:1}
An {\it $(m,n,k)$-Magog trapezoid\/} is an array of positive integers consisting
of the first $k$ rows of an array
$$
\begin{matrix} b_{11}&b_{12}&\hdotsfor2&b_{1n}\\
b_{21}&b_{22}&\dots&b_{2,n-1}\\
\hdotsfor3\\
b_{n1}
\end{matrix}
$$
such that entries along rows are weakly increasing, 
entries along columns are weakly decreasing,
and such that the entries in the first row are bounded by
$b_{11}\le m+1$, $b_{12}\le m+2$, \dots, $b_{1n}\le m+n$.
\end{definition}

\begin{figure}[h]
$$
\begin{pmatrix} 
0&0&\hphantom{-}1&\hphantom{-}0&0&0\\
1&0&-1&\hphantom{-}1&0&0\\
0&0&\hphantom{-}1&-1&0&1\\
0&1&-1&\hphantom{-}1&0&0\\
0&0&\hphantom{-}1&-1&1&0\\
0&0&\hphantom{-}0&\hphantom{-}1&0&0
\end{pmatrix}
\quad \longrightarrow\quad 
\begin{pmatrix} 
1&1&1&1&1&1\\
1&1&0&1&1&1\\
0&1&1&0&1&1\\
0&1&0&1&1&0\\
0&0&1&0&1&0\\
0&0&0&1&0&0
\end{pmatrix}
\quad \longrightarrow\quad 
\begin{matrix} 
1&2&3&4&5&6\\
1&2&4&5&6\\
2&3&5&6\\
2&4&5\\
3&5\\
4
\end{matrix}
$$
\caption{An alternating sign matrix turned into a monotone triangle}
\label{fig:Gog}
\end{figure}

Likewise, we may transform an $n\times n$ alternating sign matrix
into a triangular array. This is done by replacing the $i$-th row
of the matrix by the sum of rows $i,i+1,\dots,n$, for
$i=1,2,\dots,n$. The result is a matrix with only $0$'s and $1$'s,
see the middle of Figure~\ref{fig:Gog}. To obtain the corresponding triangular
array, for each row we record the positions of the $1$'s in that row;
see again Figure~\ref{fig:Gog}. 
The entries in the triangular array are positive integers. they
have the property that they
are strictly increasing along rows, weakly increasing 
along columns, weakly increasing along diagonals in direction
north-east,
and the entries in the first row are $1,2,\dots,n$.
Again, it is easy to see that the above defines a bijection.
The triangular arrays which one obtains here are most often called
{\it monotone triangles}.

The ``Gog trapezoids" (terminology again
stolen from \cite{ZeilBD}) below generalise these monotone
triangles. In view of the above correspondence, $n\times n$
alternating sign matrices correspond to
$(0,n,n)$-Gog trapezoids.

\begin{definition} \label{def:2}
An {\it $(m,n,k)$-Gog trapezoid\/} is an array of positive integers consisting
of the first $k$ columns of an array
$$
\begin{matrix} a_{11}&a_{12}&\hdotsfor2&a_{1n}\\
a_{21}&a_{22}&\dots&a_{2,n-1}\\
\hdotsfor3\\
a_{n1}
\end{matrix}
$$
such that entries along rows are strictly increasing, entries along
columns are weakly increasing,
and entries along 
diagonals from lower-left to upper-right are weakly increasing,
and such that the entries in the right-most column are bounded by
$a_{1k}\le m+k$, $a_{2k}\le m+k+1$, \dots, $a_{n+1-k,k}\le m+n$.
\end{definition}

Here is now the announced conjecture.

\begin{conjecture} \label{conj}
The number of $(m,n,k)$-Magog trapezoids with $s$
Maxima in the first row and $t$ Minima in the last row
equals the number of $(m,n,k)$-Gog trapezoids with $t$
Maxima in the right-most column and $s$ Minima in the left-most column. Here, a
Maximum is an entry that is equal to its upper bound, whereas a
Minimum is an entry that is $1$.
\end{conjecture}

I am convinced that there must exist a jeu-de-taquin-like 
procedure to turn the objects in Definition~\ref{def:1} 
into the ones in Definition~\ref{def:2}, 
but I am not the only one to have failed to find such a
procedure. For a different attempt to construct a bijection
between totally symmetric self-complementary plane partitions and 
alternating sign matrices (and new conjectures), see \cite{BiChAA}
and the references contained therein.

\smallskip
What is known about Conjecture~\ref{conj}? 

\smallskip
(1) If $m=0$ and $k=n$, and if we forget about Maxima and Minima,
then the conjecture reduces to the equality of the number of $n\times
n$ alternating sign matrices and the number of totally symmetric
self-complementary plane partitions contained in a $(2n)\times
(2n)\times (2n)$ box. As was said earlier, this is a theorem, thanks to
Andrews \cite{AndrAW} and Zeilberger \cite{ZeilBD}, but no
bijection is known. 

\smallskip
(2) If $m=0$ and if we forget about Maxima and Minima, then this is
Zeilberger's main theorem in \cite{ZeilBD} (namely Lemma~1), but,
again, his proof is highly non-bijective.

\smallskip
(3) If $m=0$, then the above conjecture reduces to Conjecture~7' in
\cite{MiRRAC}. 

\smallskip
(4) If $k=1$ then we want to show that the number of arrays of
positive integers
$$a_1\quad a_2\quad \dots \quad a_n$$
with $a_k\le m+k$ and with $s$ Maxima and $t$ Minima is exactly the
same as the number of arrays of the same type but with $t$ Maxima and
$s$ Minima.
It is easy to show that the number in question is
$$\binom {m+2n-s-t-2}{m+n-2}-\binom {m+2n-s-t-2}{m+n-1},$$
which is indeed symmetric in $s$ and $t$.
For this special case, it is not very difficult to construct
a bijection.

I refer the reader to \cite{RobbAA,IshiAZ,IshiBZ} for further conjectures
around totally symmetric self-complementary plane partitions
(some of them having been proved in \cite{FoZiAA}).
\medskip

In order to (finally) answer the question in the section header:
interest in plane partitions is still very much alive, if though
in different forms than at the time of MacMahon or when Stanley
studied plane partitions. Aside from the notorious problem of
connecting plane partitions and alternating sign matrices
in a conceptual way, it is
rhombus tilings and, more generally, perfect matchings of bipartite
graphs which are primarily in the focus of current research. 
On the one hand, unexpected closed form product formulae for rhombus
tilings of particular regions are constantly discovered, see e.g.\
\cite{GilmAA,LaiTAA,CoNaAA} and the references therein.
On the other hand, a deep theory of asymptotic properties of
perfect matchings (and thus, in particular, of rhombus tilings
and plane partitions) has been developed in the recent past,
see the excellent survey \cite{KenyAF}. Kasteleyn--Percus
matrices (see Section~\ref{sec:methods}) are the starting point
of this theory. To go into more detail
of the subsequent asymptotic analysis
would however definitely go beyond the scope of these notes, 
and therefore this is the place to stop.

\end{document}
